%% file: character-final.tex
\documentclass[11pt]{amsart}

\usepackage{url}
\usepackage[utf8]{inputenc}
\usepackage{asymptote}
\usepackage{latexsym,amssymb,amsmath}
\usepackage{mathrsfs}
\usepackage{epsfig}
\usepackage{amsmath}
\usepackage{amsfonts}
\usepackage{amssymb}
\usepackage{amssymb,amscd}
\usepackage{graphicx}
\usepackage{subfigure}
\usepackage{psfrag}
\usepackage{comment}
\usepackage[all]{xy}\xyoption{rotate}\xyoption{pdf}

\RequirePackage{graphics}

\usepackage{bbm}
\usepackage[english]{babel}

\usepackage[left,modulo]{lineno}
\catcode`\<=\active \def<{
\fontencoding{T1}\selectfont\symbol{60}\fontencoding{\encodingdefault}}
\catcode`\>=\active \def>{
\fontencoding{T1}\selectfont\symbol{62}\fontencoding{\encodingdefault}}
\catcode`\|=\active \def|{
\fontencoding{T1}\selectfont\symbol{124}\fontencoding{\encodingdefault}}

\newcommand{\nosymbol}{}
\newcommand{\tmem}[1]{{\em #1\/}}
\newcommand{\tmop}[1]{\ensuremath{\operatorname{#1}}}
\newcommand{\tmtextbf}[1]{{\bfseries{#1}}}

\newtheorem{theorem}{Theorem}

\newtheorem{prop}{Proposition}
\newtheorem{lem}[theorem]{Lemma}

\newtheorem{pro}{Proposition}[section]
\theoremstyle{definition}
\newtheorem{defn}{Definition}[section]

\newcommand\bC{\mathbb{C}}

\newcommand\SL{\textrm{SL}}
\newcommand\PSL{\textrm{PSL}}

\newcommand\PGL{\textrm{PGL}}
\newcommand\Geight{\Gamma_{8}}
\newcommand{\Hom}{\mathrm{Hom}}

\newcommand{\ind}{\setlength{\parindent}{2em}}
\newcommand{\noind}{\setlength{\parindent}{0em}}

\def\bfi#1#2{\fontsize{#1}{#2}\fontseries{b}\fontshape{it}\selectfont}
\def\bc{\,,\,}
\def\h{\hspace{.5em}}

\def\vv{\vspace{1ex}}
\def\vvv{\vspace{2ex}}

\def\square{\hfill${\vcenter{\vbox{\hrule height.4pt \hbox{\vrule width.4pt
height7pt \kern7pt \vrule width.4pt} \hrule height.4pt}}}$}

\newenvironment{pf}{{\it Proof.}\quad}{\square \vskip 12pt}

\author{E. Falbel, A. Guilloux, P.-V. Koseleff,\\ F. Rouillier and M. Thistlethwaite}

\begin{document}

\title[] {Character varieties for $\SL(3,\bC)$: the figure eight knot. }

\begin{abstract}
We give a description of several representation varieties of the fundamental group of the complement of the figure eight knot in
$\PGL(3,\bC)$ or $\PSL(3,\bC)$.  We  obtain a description of the projection of the representation variety into
the character variety of the boundary torus into $\SL(3,\bC)$. 
\end{abstract}

\maketitle
\input{character-11}
\input{character-2new.tex}
\input{character-3_mbt}

\end{document}

%% file: character-11.tex
\section*{Introduction}

Representation varieties of finitely presented  groups into Lie groups has been studied since a long time.  A recent by Sikora paper with many references
is  \cite{Sikora}.  In this paper we concentrate on the representations of 3-manifold groups in $\PGL(3,\bC)$ or $\SL(3,\bC)$.
  The case of surface groups have been studied more extensively and their  representation varieties
into $\SL(3,\bC)$ are treated by Lawton in \cite{Lawton}.  On the other hand, only recently the study
of representations of 3-manifold groups into $\SL(3,\bC)$ or  $\PGL(3,\bC)$  has been started. See for instance the work of Porti and Menal-Ferer \cite{MF}, of Garoufalidis, Goerner, Thurston and Zickert \cite{GTZ,GGZ} and of Bergeron with the two first named authors \cite{BFG} among others.  One has to keep in mind
the deep study which was started in the last decades of the last century on representations into $\SL(2,\bC)$ by several authors among them
Thurston \cite{Thurston}.  It gives a paradigm for the present study in the case of  $\PGL(3,\bC)$.  

Let $\Geight$ be the fundamental group of the figure eight knot complement. We present here a thorough study of the space of representations of $\Geight$ into $\PGL(3,\bC)$ and $\SL(3,\bC)$.  It turns out that all the irreducible representations we find in  $\PGL(3,\bC)$ can be lifted to $\SL(3,\bC)$.  Therefore, most of the computations will be made in $\PGL(3,\bC)$ but we will occasionally work with $\SL(3,\bC)$.

We begin with a review of different versions of this space (representation variety, character variety, decorated representations among others) and the coordinates we use to study these varieties in a rather general setting.  A diagram, displayed in Figure \ref{fig:diagram}, shows the several representation spaces and their relations.  The main coordinate system we use is derived from projective flag decorations of ideal triangulations of a 3-manifold as introduced in \cite{BFG} (see also \cite{GGZ}).   Another view of these spaces is given in  \cite{Zickert-A-pol} using Ptolemy coordinates, which are coordinates for affine flags decorations of ideal triangulations. Those work are accessible through the website CURVE \cite{CURVE}.

Then we proceed in section \ref{computations} to an explicit description of these different spaces for the group $\Geight$ and the different maps relating these spaces.  We use coordinates on a Zariski open set of our spaces which are associated to a given triangulation.  In particular, we describe 
a Zariski open set (called the deformation variety) of the set of decorated representations into $\PGL(3,\bC)$ given by coordinates we introduced in \cite{BFG}.
The main theorem is  Theorem \ref{3components}: 
Given the standard triangulation of the figure eight knot  complement, there exits precisely three irreducible components of the
 deformation variety.  Each component is smooth of dimension two. Whereas this result is stated and proved only for the deformation variety, we strongly believe this is also the case for the whole character variety. A precise statement and proof should appear elsewhere.
This section begins with a reminder about Groebner basis and saturation, then a presentation of the actual equations for the figure eight knot complement in \S \ref{ss:equations}. We then proceed with a general presentation of our method to achieve computations. Eventually, we give precisions of the specific case of the figure eight knot complement.
The three components are finally described in subsection \ref{results}.  It turns out that each representation can be lifted to an $\SL(3,\bC)$ representation. Once these components are identified, we may check the properties stated in the theorem (irreducibility, smoothness...).

The next section \ref{s:A-var} is devoted to some more precise understanding of these components: description the $A$-variety associated to this knot (it is a natural analog of the A-polynomial), identification of the boundary unipotent representation and identification of the space of representations in $\PGL(2,\bC)$.

Finally, in section \ref{deformation}, we describe a deformation method to obtain a parametrization of some generators of $\Geight$ on each irreducible components.
This method allows one to obtain a particular irreducible component
of the representation variety containing a given representation.  It is based on a Newton method and LLL algorithm and was successfully used in other contexts
\cite{CLTa,CLTb}. We applied it to the boundary unipotent representations
of the fundamental group of the figure eight knot complement. A great interest of this method is that its output is indeed a parametrization of matrices. It opens the path to a geometric study of those representations, see e.g. \cite{DF}.
However it does not belong to its scope to determine all the irreducible components. We insist on the fact that it is the conjunction of both methods that gives a satisfactory result: a parametrization of each component.

We thank Nicolas Bergeron, Martin Deraux,  Matthias Goerner, Julien March\'e, Pierre Will, Maxime Wolff and Christian Zickert for many discussions and exchanges during this project.
The authors are grateful to the ANR 'Structures G\'eom\'etriques Triangul\'ees' which has financed this project during  the last two years.

\section{Representation spaces}

We begin with a review of the different spaces classically considered as representation spaces. As we will consider several variations of the same space, we try to use explicit -- but sometimes lenghty -- names for these spaces.

Let $\Gamma$ be a finitely presented group and $G$ the group of points over $\bC$ of a linear algebraic group. In the following, we will mainly consider $G=\PGL(3,\bC)$ or $\SL(3,\bC)$, but we will occasionally look also at $\PGL(2,\bC)$, $\SL(2,\bC)$ or even the real group $\mathrm{PU}(2,1)$. 

\subsection{Representation variety and character variety}

We begin by the most classical representation variety and character variety.

\newcommand{\Mon}{\mathrm{Mon}}

\begin{defn} 
We denote by $\Hom(\Gamma, G)$ the set of all morphisms
from $\Gamma$ to $G$. This space is called the \emph{representation variety}.

The \emph{monodromy variety}\footnote{Neither the name nor the notation is classical. The name has been so choosen because a geometric structure on a variety gives naturally, via the monodromy of the structure, an element in the monodromy variety.} $\Mon(\Gamma,G)$ is the quotient under the natural action by conjugation of $G$ :
$$\Mon(\Gamma,G)=\Hom(\Gamma,G)/G.$$
\end{defn}

Note that $\Hom(\Gamma,G)$ is an affine algebraic variety as $\Gamma$ is finitely presented and $G$ is affine  algebraic. Moreover the relators correspond to generators of the ideal defining the variety. Note however that $\Mon(\Gamma,G)$ is not an algebraic variety.

Though, the action of $G$ on itself by conjugation is an algebraic action.  Therefore it defines an algebraic action on $\Hom(\Gamma, G)$ which induces an
algebraic action on its regular functions.  We let $\bC[\Hom(\Gamma, G))]^{G}$ be the ring of invariant functions.
\begin{defn}
The \emph{character variety} is the algebraic quotient:
$$
X(  \Gamma, G)=\Hom(\Gamma, G)//G.
$$
It is the affine variety associated to the ring  $\bC[\Hom(\Gamma, G)]^{G}$ and comes together with the  regular map
$\Hom(\Gamma, G)\rightarrow  X(\Gamma, G)$ induced by the homomorphism 
$$
  \bC[X(  \Gamma, G)]=\bC[\Hom(\Gamma, G)]^{G} \rightarrow \bC[\Hom(\Gamma, G)]. 
$$
\end{defn}
\newcommand{\tr}{\mathrm{tr}}
The name of this variety comes from its links with the set of \emph{characters}, at least when $G=\SL(n,\bC)$ or $\PGL(n,\bC)$. 

When $G=\SL(n,\bC)$ traces of elements are well defined and one knows that there is a bijection between the character variety and the set of trace functions.

For the case $G=\PGL(n,\bC)$, let us remark that for an element $g$ in $\PGL(n,\bC)$, one may define its trace to the power $n$: choose a representant $\tilde g$ of $g$ which belongs to $\SL(n,\bC)$. Then we put: $\tr^n(g):=\tr^n(\tilde g)$. It does not depend on the choice of $\tilde g$.

\begin{defn} The character of $\rho\in\Hom(\Gamma,G)$ is the trace function
$\chi_\rho : \Gamma\rightarrow \bC$ defined by $\chi_\rho(\gamma)={\mbox{tr}}^n(\rho(\gamma))\ $. 
\end{defn}
 
Although we will not use it here, we believe that the map 
$$ 
X(  \Gamma, G)\rightarrow \{\textrm{characters of }\Gamma\} 
$$
induced by $\rho \mapsto \chi_\rho$ is a bijection. It is well-known for $n=2,\, 3$, see e.g. \cite{Lawton}.
Note that there is also a natural application from $\Mon(\Gamma,G)$ to $X(\Gamma,G)$ sending the class of a representation $\rho$ to $\chi_\rho$. This application is a surjection, but not an injection. It is not algebraic, as $\mathrm{Mon}(\Gamma,G)$ is not an algebraic variety.

\subsection{Decorated versions}

From now on, we assume $G=\PGL(n,\bC)$. A key point to study the previous defined spaces is to add a geometric structure on it, as usual in the study of moduli spaces. 
Such a structure is called a \emph{decoration} and will be most easily defined with additional assumption: from now on $\Gamma$ is the fundamental group of a $1$-cusped hyperbolic $3$-variety $M$ (as the figure eight knot complement)\footnote{This assumption is not necessary and one may work with a triangulated $3$-variety (either hyperbolic with more cusps, or with boundary, or even not hyperbolic). We refer to \cite{BFG} for this setting.}. As $M=\mathbb H^3/\Gamma$ is a cusped hyperbolic manifold, we let
 $\mathcal P\subset \mathbb{CP}^1=\partial \mathbb H^3$ be the set of cusps in the boundary of the hyperbolic space. By construction, $\Gamma$ acts on $\mathcal P$.
Let moreover $\mathcal Fl$ be the space of complete flags of $\mathbb{CP}^n$. There is a natural action of $\PGL(n,\bC)$, for which $\mathcal Fl$ becomes the quotient of $\PGL(n,\bC)$ by the subgroup of upper-triangular matrices. Hence, given a representation $\rho$ of $\Gamma$ in $\PGL(n,\bC)$, we have an action of $\Gamma$ on $\mathcal Fl$, through $\rho$.  For these objects and the following, we refer to \cite{BFG}.

\newcommand{\DecHom}{\mathrm{DecHom}}

\begin{defn}
Let $\rho \in \Hom(\Gamma,G)$. Then $\Gamma$ acts both on $\mathbb{CP}^1$ and $\mathcal Fl$. A \emph{decoration} of $\rho$ is a $\Gamma$-covariant map :
$$\phi : \mathcal P\to \mathcal Fl.$$
Let $\DecHom(\Gamma,G)$ be the decorated representations variety:
$$\DecHom(\Gamma,G)=\{(\rho,\phi)|\; \rho\in \Hom(\Gamma,G), \; \phi\textrm{ is a decoration of }\rho\}.$$ 
\end{defn}

As the group $G=\PGL(n,\bC)$ acts on $\mathcal Fl$, it acts on the space of maps $\mathcal P\to \mathcal Fl.$
Moreover, it is easily checked that if $\phi$ is a decoration of $\rho$, then $g\cdot \phi$ is a decoration of the conjugated $g\cdot \rho=g\rho g^{-1} $. 
Hence $G$ acts on the space $\DecHom(\Gamma,G)$. 

\newcommand{\DecMon}{\mathrm{DecMon}}

\begin{defn}
We define the \emph{decorated monodromy varieties} $\DecMon(\Gamma,G)$ as the naive quotient:
$$\DecMon(\Gamma,G)=\DecHom(\Gamma,G)/G.$$
\end{defn}

We now proceed by defining a decorated version of the character variety. 
First of all, note that the space of flags $\mathcal Fl$ is an algebraic variety and the action of $G$ is algebraic. 
Moreover, $\DecHom(\Gamma,G)$ is an algebraic variety: indeed, choose a \emph{finite} fundamental set $\{p_1,\ldots,p_r\}$ for the action of $\Gamma$ on $\mathcal P$. 
Then a decorated representation $(\rho,\phi)$ is uniquely determined by a choice of a flag $F_i$ for each $p_i$ such that $F_i$ is stabilized by $\mathrm{Stab}_\Gamma(p_i)$.

\newcommand{\DecX}{\mathrm{Dec}X}

We now consider the algebraic quotient of $\DecHom(\Gamma,G)$ under this action:
\begin{defn}
We define the \emph{decorated character variety} $\DecX(\Gamma,G)$ as the algebraic quotient:
$$\DecX(\Gamma,G)=\DecHom(\Gamma,G)//G.$$
\end{defn}

\subsection{Peripheral representations: choosing a meridian and a longitude}\label{ss:lm}
\newcommand{\Hol}{\mathrm{Hol}}\newcommand{\periph}{\mathrm{periph}}

By considering decorated versions of the representations, we also grant an easy parametrization 
for the \emph{peripheral representation}, that is the restriction of the representation to the fundamental group of the peripheral torus. This is more carefully studied in \cite[Section 1]{G}. Recall that $M$ is a $1$-cusped hyperbolic $3$-manifold. Let $\mathbb T$ be the peripheral torus.

First, let $\Gamma_\mathrm{para}$ be the set of parabolic elements of $\Gamma$ (the identity is not considered parabolic). 
Each element $\gamma$ of $\Gamma_\mathrm{para}$ stabilizes a unique cusp $p_\gamma\in\mathcal P$. 
Consider a decorated representation $(\rho,\phi)\in \DecHom(\Gamma,G)$. 
Let $D$ be the diagonal group of $G$. Then, in any basis adapted to the flag $\phi(p_\gamma)$,
$\rho(\gamma)$ is an upper-diagonal matrix. Moreover, the diagonal part of $\rho(\gamma)$ is invariant under the action of $G$ \cite{G}.
Hence we get a $G$-invariant map: 
$$\DecHom(\Gamma,G)  \to \{\Gamma_\mathrm{para}\to D\}.$$
\newcommand{\dec}{\mathrm{dec}}

In order to parametrize the peripheral representations, we fix, \emph{once for all}, the following choices: 
\begin{itemize}
\item $\pi_1(\mathbb T)\simeq \mathbb Z^2$, through the choice of a longitude $l$ and a meridian $m$ of $\mathbb T$.
\item An injection $\pi_1(\mathbb T)\to \Gamma=\pi_1(M)$.
\end{itemize}

A representation in $\Hom(\Gamma,G)$ may be restricted to $\pi_1(\mathbb T)$. This gives an algebraic map of \emph{restriction}:\newcommand{\restr}{\mathrm{restr}}
$$X(\Gamma,G)\xrightarrow{\restr}X(\mathbb Z^2,G).$$
Moreover, in the decorated case, as $\pi_1(\mathbb T)$ is composed of parabolic elements, we get a map, called \emph{peripheral holonomy map} and denoted by $\Hol_\periph$:
$$\DecX(\Gamma,G)\xrightarrow{\Hol_\periph}X(\mathbb Z^2,D)=:\DecX(\mathbb Z^2,G).$$
It consists, at the level of decorated representations, in sending a decorated representation to the diagonal part of the restriction restricted to $\pi_1(\mathbb T)$.

The last step of the parametrization is achieved by noting that $D$ is an affine algebraic group isomorphic to $(\bC^*)^{n-1}$ in the case $G=\PGL(n,\bC)$ or $\SL(n,\bC)$. We fix \emph{once for all} such an isomorphism. Using this, the space $\DecX(\mathbb Z^2,G)=X(\mathbb Z^2,D)$ is isomorphic to $(\bC^*)^{2(n-1)}$. This is the desired parametrization of the peripheral representations.

 In the case of $\PGL(3,\bC)$, we choose (for computational reasons, see section \ref{computations}) more precisely to use the isomorphism between $D$ and $(\bC^*)^2$ given by:
$$\begin{pmatrix} L^*&&\\&1&\\&&L\end{pmatrix}\quad \mapsto\quad (L,L^*).$$
 So the space $X(\mathbb Z^2,D)$ is isomorphic to $(\bC^*)^4$. A representation such that the diagonal part of the longitude and the meridian are:
 \[
\rho(l) = \left[ \begin{matrix}
  L^* & \ast & \ast\\
    0 & 1 & \ast\\
    0 & 0 & L
  \end{matrix} \right]\quad \textrm{and}\quad \rho(m) = \left[ \begin{matrix}
  M^* & \ast & \ast\\
    0 & 1 & \ast\\
    0 & 0 & M
  \end{matrix} \right]
 \]
will be represented by the coordinates $(L,L^*,M,M^*)$ (the choice is made to be compatible with the one 
we used in \cite{BFG}).

\subsection{Deformation variety}\label{ss:deformationvariety}

Thurston \cite{Thurston} defined coordinates on 
$$\DecX(\Gamma,\PGL(2,\bC))
$$ by defining the \emph{deformation variety} with the \emph{gluing equations}. We use in this paper an analog og this space for the group $\PGL(3,\bC)$. From now on, we will always assume that $n=3$. This space is defined by further assuming that $M$ is ideally triangulated: $M$ is homeomorphic to a gluing of tetrahedra with vertices removed as in \cite{Thurston, BFGKR}
So let $\mathcal T=(T_1,T_2,\ldots,T_\nu)$ be an ideal triangulation of $M$.

\begin{figure}
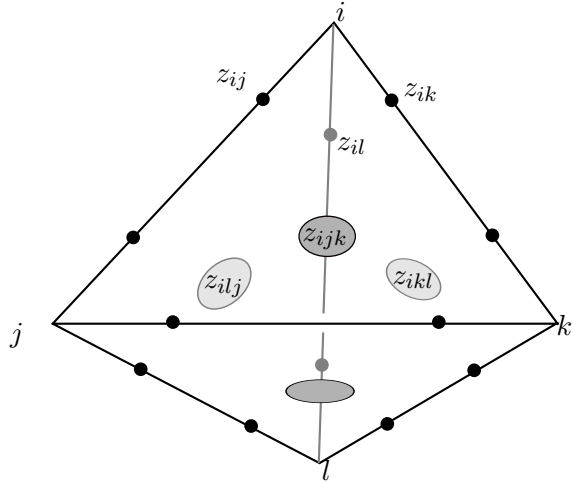
\caption{Coordinates for a tetrahedron of flags}\label{coordinates}
\end{figure}

We now recall coordinates on these spaces as defined in \cite{BFG} (see Figure \ref{coordinates}).  For each of these (oriented) tetrahedra, one consider a set of $16$ coordinates: one on each half-edge and one in the center of this face. If the vertices of the tetrahedron $T$ are $(ijkl)$, one has twelve coordinates on half edges, denoted by $z_{ij}(T)$, $z_{ik}(T)$... and four on faces, denoted by $z_{ijk}(T)$.

 As described in \cite{BFG}, these coordinates are subject to two different kinds of consistency relations: 
\begin{itemize}
 \item Internal relations: for a tetrahedron $T=(ijkl)$ we have around each vertex $z_{ij}(T)z_{ik}(T)z_{il}(T)=-1$ and $z_{ik}(T)=\frac{1}{1-z_{ij}(T)}$; and for each face $z_{ijk}(T)=-z_{il}(T)z_{kl}(T)z_{jl}(T)$.  
 \item Gluing relations: first, given two adjacent tetrahedra $T_{\mu}$, $T_{\mu'}$
of $T$ with a common face  
$(ijk)$ then 
\begin{equation} \label{face}
z_{ijk} (T_{\mu}) z_{ikj} (T_{\mu'}) = 1.
\end{equation}
second given a sequence $T_1, \ldots , T_{\mu}$ of tetrahedra sharing a
common edge $ij$ and such that $ij$ is an inner edge of the sub-complex composed by $T_1\cup \cdots \cup T_{\mu}$ 
then 
\begin{equation} \label{edge}
z_{ij} (T_1) \cdots z_{ij} (T_{\mu}) = z_{ji} (T_1) \cdots z_{ji} (T_{\mu}) = 1.
\end{equation}
 \end{itemize}
 
Given $M$ and a triangulation $\mathcal{T}$ with $n$ tetrahedra, we consider the space of solutions of these equations and denote it by 
$$
\mathrm{Defor}(M,\mathcal{T})\subset \bC^{16n}.
$$
We call it the {\it deformation variety} of a triangulation. 

Up to combinatoric choices, the holonomy map associates to each point in $\mathrm{Defor}(M,\mathcal T)$ a decorated representation into $\PGL(3,\bC)$ as in \cite[section 5]{BFG}: 
$$
\mathrm{Defor}(M,\mathcal{T})\rightarrow \mathrm{DecHom}(\pi_1(M),\PGL(3,\bC)).
$$
Moreover this map is algebraic\footnote{indeed, by construction, the expression of the matrix entries are algebraic.} and a different initial choice gives a conjugated decorated representation.  

So this gives two well defined maps, landing in $\DecMon$ or $\DecX$, both called holonomy maps and denoted by $\Hol$. Note that the second one is algebraic: 
$$
\Hol\: : \:\mathrm{Defor}(M,\mathcal{T})\rightarrow \DecMon(\pi_1(M),\PGL(3,\bC)),
$$
and
$$
\Hol\: : \:\mathrm{Defor}(M,\mathcal{T})\rightarrow \DecX(\pi_1(M),\PGL(3,\bC)),
$$
We extend the peripheral holonomy map defined in the previous section in a natural way:
\[
\Hol_\periph\: : \:\mathrm{Defor}(M,\mathcal{T})\rightarrow\DecX(\pi_1(T),\PGL(3,\bC))\simeq(\bC^*)^4
\]

Let us comment a bit on the properties of both the holonomy maps $\Hol$: the first map is injective but not surjective, as some representation are not detected by a given triangulation. Moreover, if you forget the decoration, the maps are not any more injective : there is a (generically) finite choice for a decoration of a given representation.  All this can be seen already in the case of representations
into $\PGL(2,\bC)$. In particular the identity representation can be obtained by gluing two tetrahedra to obtain the 
sphere $S^3$ minus  four points and the solutions to the hyperbolic equations has a curve such that the associated representations are 
trivial \cite{SegTill}.  Also, in the tables obtained in \cite{BFGKR} the 1-dimensional components give conjugated representations with  finite group image into $\PGL(3,\bC)$.
Examples where a given triangulation is not enough to describe all representation are given in \cite{SegTill} already in the case of
representations into $\PGL(2,\bC)$.

On the other hand, the situation is good on some components. First, recall that, as $M$ has been assumed to be hyperbolic, there is the monodromy $[\rho_{\mathrm{geom}}]\in\Mon(\Gamma,\PGL(2,\bC))$ of the unique oriented complete hyperbolic structure on $M$. We call it the \emph{geometric representation}. Moreover, the irreducible representation $\PGL(2,\bC)\to\PGL(n,\bC)$ gives a map $\Mon(\Gamma,\PGL(2,\bC))\to\Mon(\Gamma,\PGL(n,\bC))$. We still call the image of $[\rho_{\mathrm{geom}}]$ the \emph{geometric representation}, still denoted by $[\rho_{\mathrm{geom}}]\in\Mon(\Gamma,\PGL(n,\bC))$. One can then state the following theorem for $\PGL(2,\bC)$ see \cite{Thurston,SegTill}, with a mild assumption on the triangulation that we will not define here and that is verified for the triangulation of the figure eight knot complement we will use:

\begin{theorem} Let M be an orientable, connected, cusped hyperbolic manifold.  Let $\mathcal{T}$
be an ideal triangulation of M such that all edges are essential.  Then there exists $z_{\mathrm{geom}}\in \mathrm{Defor}(M,\mathcal{T})$
such that $\Hol(z_{\mathrm{geom}})$ is the geometric representation.  

Moreover the whole  component of $X(\Gamma,\PGL(2,\bC))$ containing the (image of) the geometric representation is in the image  of the holonomy map $\Hol$.
\end{theorem}

The analogous theorem in the case of $\PGL(3,\bC)$ is believed to be true:  
Let $M$ be an orientable, connected, cusped hyperbolic manifold.  Let $\mathcal{T}$
be an ideal triangulation of $M$ such that all edges are essential.  Then there exists $z_{\mathrm{geom}}\in \mathrm{Defor}(M,\mathcal{T})$
such that $\Hol(z_{\mathrm{geom}})$ is the geometric representation. We believe that the image of $\Hol$ contains the whole  component of $X(\Gamma,\PGL(3,\bC))$ containing the geometric representation.

We also know that $z_{\mathrm{geom}}\in \mathrm{Defor}(M,\mathcal{T})$ is a smooth point (see \cite{BFGKR}).

\subsection{A diagram}

At this point, one may write a commutative diagram: see figure \ref{fig:diagram}. Note that, between algebraic varieties, the maps are algebraic. Moreover, the vertical maps are just the forgetful maps. The non-algebraic maps are denoted by dashed lines. For the sake of readability, we drop the mentions of $\Gamma$ and $G$. Moreover, we denote by $\DecX_{\mathbb T}$ and $X_{\mathbb T}$ the (decorated) character variety of $\pi_1(\mathbb T)$. Recall that we have done a choice of an isomorphism $D\simeq (\bC^*)^2$ and of a longitude $l$ and meridian $m$ in the torus $\mathbb T$, so that $\DecX_{\mathbb T}$ is isomorphic to $(\bC^*)^4$.
\begin{figure}[ht]
\begin{center}
\xymatrixcolsep{3.5pc}\xymatrix{
\mathrm{Defor} \ar@/^2.5pc/[ddrrr]^{\: \Hol \:} \ar[ddr]_{\small \txt{(up to\\choices)}} \ar@{-->}[drr]^{\: \Hol\: } \ar@/^3pc/[ddrrrr]^{\: \Hol_\periph \:}&&&&&\\
&&\DecMon \ar@{-->}[dr] \ar@{-->}'[d]'[dd][ddd]&&&\\
&\DecHom \ar[d] \ar@{-->}[ur] \ar[rr]&&\DecX \ar[d] \ar[r]^-{\Hol_\periph}&\DecX_{\mathbb T}\ar[d]\\
&\Hom \ar[rr] \ar@{-->}[dr] &&X \ar[r]^{\restr}&X_{\mathbb T}\\
&&\Mon \ar@{-->}[ur]&&&
}
\end{center}
\caption{The different representation varieties}\label{fig:diagram}
\end{figure}

\subsection{The A-variety}

In this section we will work with representations with values in $\SL(3,\bC)$. As before $M$ is a 3-manifold with boundary a torus where we fix a basis 
of the homology group given by a choice of a longitude and a meridian. One can identify $\DecX(\mathbb Z^2,G)$ to the diagonal representations of $\mathbb Z^2$.  As before, we see that
$\DecX(\mathbb Z^2,G)$ is isomorphic to $(\bC^*)^4$.  The map $\DecX_{\mathbb T}\rightarrow X_{\mathbb T}$ is the quotient by the permutation group $S_3$ acting on the triple of complex numbers whose product is $1$.

Consider the closure of the image of $\DecX(\pi_1(M),\SL(3,\bC))$ in $\DecX(\mathbb Z^2,G)$.   Let $D_M$  be the union of the component of maximal dimension of this closure.

Consider now a natural embedding  ${\bC^*}^4\subset \bC P^4$ given by ${(L,L^*,M,M^*,1)}$ and $\bar D_M$ the closure of the 
image of $D_M$. The ideal boundary of $D_M$ is $\bar D_M\setminus D_M$.  Essentially the following definition was given independently in \cite{Zickert-A-pol}. 

\begin{defn}\label{defaideal}
We define the $A$-variety of $M$ (for $n=3$ and with a choice of basis of the boundary torus homology) to be the closure of $D_M$ in  $\bC^4$.
We define the \emph{A-ideal} to be its defining ideal.
\end{defn}

\section{The figure eight knot complement}\label{section:intro-eight}

We present here the geometric and combinatorial facts on the figure eight knot complement $M$ we use afterwards. We need a triangulation $\mathcal T$ of $M$ -- in order to study the variety $\mathrm{Defor}(M,\mathcal T)$ -- together with some presentations of its fundamental group for paramtrizing matrices.

\subsection{Triangulation}

\begin{figure}[ht]
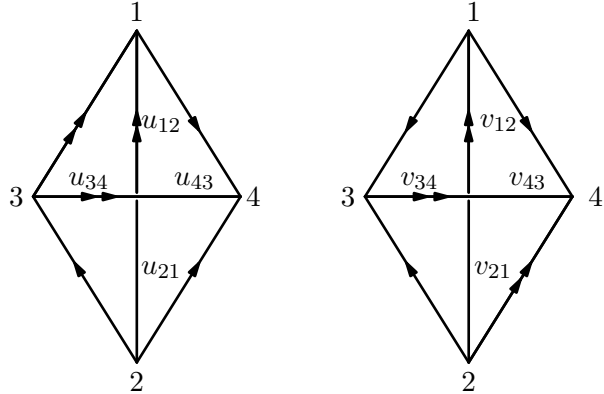
 \label{eight}
\begin{center}
\begin{asy}
   size(8cm);
  defaultpen(1);
  usepackage("amssymb");
  import geometry;
  
  point o = (8,0);
  point oo = (24,0);
  pair ph = (10,0);
  pair pv = (5,8);

  // 
  draw(o -- o+pv,Arrow(5bp,position=.4), Arrow(5bp,position=.5));
  draw(o+pv -- o+ph, Arrow(5bp,position=.6));
  draw(o+ph-pv-- o, Arrow(5bp,position=.6));
  draw(o+ph-pv-- o+ph, Arrow(5bp,position=.6));
  draw(o+ ph/2+(0,.2)--o+pv, Arrow(5bp,position=.4), Arrow(5bp,position=.5));
  draw(o+ph-pv -- o+ph/2+ (0,-.2));
  draw(o -- o+ph,Arrow(5bp,position=.3), Arrow(5bp,position=.4) );
  
  // 
  label("\small { $u_{34}$}", o, 3*dir(25));
  //
  //
  label("{\small  $3$}", o, 1*W);
  
  label("\small {$u_{43}$}", o+ph, 3*dir(180-25));
  //
  //
  label("{\small $4$}", o+ph, .6*E);
  
  //
  label(" \small $ u_{12}$", o+pv, 9*dir(-80));
  //
  label("{\small $1$}", o+pv, 1*N);
  
   //
  label(" \small $ u_{21}$", o+ph-pv, 9*dir(80));
  //
  label("{\small $2$}", o+ph -pv, 1*S);

  // 
  draw( oo+pv -- oo,Arrow(5bp,position=.6));
  draw( oo+pv -- oo+ph,Arrow(5bp,position=.6));
  draw( oo+ph-pv -- oo,Arrow(5bp,position=.6));
  draw(oo+ ph/2+(0,.2)--oo +pv,Arrow(5bp,position=.4), Arrow(5bp,position=.5));
  draw(oo+ph-pv -- oo+ ph/2-(0,.2));
  draw(oo+ph-pv -- oo+ ph,Arrow(5bp,position=.5),Arrow(5bp,position=.6));
  draw(oo -- oo+ph,Arrow(5bp,position=.3), Arrow(5bp,position=.4));
  
   // 
  label("\small { $v_{34}$}", oo, 3*dir(25));
  //
  //
  label("{ \small $3$}", oo, 1*W);
  
  label("\small {$v_{43}$}", oo+ph, 3*dir(180-25));
  //
  //
  label("{ \small $4$}", oo+ph, .6*E);
  
  //
  label(" \small $  v_{12}$", oo+pv, 9*dir(-77));
  //
  label("{ \small $1$}", oo+pv, 1*N);
  
   //
  label(" \small $ v_{21}$", oo+ph-pv, 9*dir(80));
  //
  label("{\small $2$}", oo+ph -pv, 1*S);

\end{asy}
\caption{ The figure eight knot represented by two tetrahedra.} \label{fig:triang8}
\end{center}
\end{figure}

It is a well-known fact, due to Riley and used by Thurston \cite{Thurston}, that the figure eight knot complement may be triangulated by two tetrahedra, with the combinatorics of face gluings displayed in figure \ref{fig:triang8}.
In order to simplify the notations, we denote by $u_{ij}$ and $v_{ij}$ be the coordinates associated to the edge $ij$ of the two tetrahedra as  shown in Figure \ref{fig:triang8}.

From this triangulation, it is an elementary exercise to compute the equations defining $\mathrm{Defor}(M,\mathcal T)$ as explained in \S \ref{ss:deformationvariety}. These equations will actually be displayed in \S \ref{ss:equations}.

\subsection{Presentations}

Depending on the situations we will use two different presentations of the figure eight knot complement fundamental group $\Gamma_8$. 

We first use the presentation (with generators as in \cite{FKR14})
$$ 
\langle\ g_1, g_2,
g_3\ |\ g_2^{-1}g_1g_3^{-1}g_1^{-1}g_3,\ g_3^{-1}g_2g_1g_2^{-1}\ \rangle.
$$  
Observe then that
$$
g_2 = [g_1,g_3^{-1}]
\quad
\textrm{and}\quad
g_3 = g_2g_1g_2^{-1}.
$$
This presentation is the SnapPea non-simplified presentation with a change of notations.
Keeping only the generators $g_1$ and $g_3$, we get the usual parabolic presentation, used in \S \ref{ss:R_2-R_3}, with $a=g_1$ and $b=g_3$.

The simplified presentation in SnapPea, which we use in \S \ref{ss:R1}, is given (uppercase denotes inverse) by
$$ 
\langle\ a, b\ |\ \ ab^3aBA^{2}B \ \rangle
$$  
with $a=  g_3g_2$ and $b=g_2^{-1}$ (hence $g_1=ba$).

With this notation, we know \cite{Thurston} that the canonical meridian is $m=ab=g_3$ and the longitude 
$l=abABaBAb=[g_3,g_2^{-1}][g_3,g_2]$.




%% file: coordinates.pdf_tex
\begingroup%
  \makeatletter%
  \providecommand\color[2][]{%
    \errmessage{(Inkscape) Color is used for the text in Inkscape, but the package 'color.sty' is not loaded}%
    \renewcommand\color[2][]{}%
  }%
  \providecommand\transparent[1]{%
    \errmessage{(Inkscape) Transparency is used (non-zero) for the text in Inkscape, but the package 'transparent.sty' is not loaded}%
    \renewcommand\transparent[1]{}%
  }%
  \providecommand\rotatebox[2]{#2}%
  \ifx\svgwidth\undefined%
    \setlength{\unitlength}{222.04707031bp}%
    \ifx\svgscale\undefined%
      \relax%
    \else%
      \setlength{\unitlength}{\unitlength * \real{\svgscale}}%
    \fi%
  \else%
    \setlength{\unitlength}{\svgwidth}%
  \fi%
  \global\let\svgwidth\undefined%
  \global\let\svgscale\undefined%
  \makeatother%
  \begin{picture}(1,0.81858688)%
    \put(0,0){\includegraphics[width=\unitlength]{coordinates.pdf}}%
    \put(0.54895741,0.78573912){\color[rgb]{0,0,0}\makebox(0,0)[lb]{\smash{$i$}}}%
    \put(-0.00358877,0.24099334){\color[rgb]{0,0,0}\makebox(0,0)[lb]{\smash{$j$}}}%
    \put(0.92347484,0.25205013){\color[rgb]{0,0,0}\makebox(0,0)[lb]{\smash{$k$}}}%
    \put(0.52538854,0.00635423){\color[rgb]{0,0,0}\makebox(0,0)[lb]{\smash{$l$}}}%
    \put(0.34515013,0.68201553){\color[rgb]{0,0,0}\makebox(0,0)[lb]{\smash{$z_{ij}$}}}%
    \put(0.55077708,0.56737403){\color[rgb]{0,0,0}\makebox(0,0)[lb]{\smash{$z_{il}$}}}%
    \put(0.66359914,0.65835936){\color[rgb]{0,0,0}\makebox(0,0)[lb]{\smash{$z_{ik}$}}}%
    \put(0.49254638,0.41269876){\color[rgb]{0,0,0}\makebox(0,0)[lb]{\smash{$z_{ijk}$}}}%
    \put(0.3251333,0.33263156){\color[rgb]{0,0,0}\makebox(0,0)[lb]{\smash{$z_{ilj}$}}}%
    \put(0.64482958,0.33950299){\color[rgb]{0,0,0}\makebox(0,0)[lb]{\smash{$z_{ikl}$}}}%
  \end{picture}%
\endgroup%

%% file: character-2new.tex
\section{The deformation variety for the figure eight knot complement}
\label{computations}

We fix the usual ideal $2$-tetrahedra triangulation of the figure eight knot
complement (which we refer as the standard triangulation). In this section we
prove the following theorem.

\begin{theorem}\label{3components}
  Given the standard triangulation of the figure eight knot complement, the
  deformation variety $\tmop{Defor} ( M, \mathcal{T} )$ is the union of 3
  distinct smooth affine (irreducible algebraic) varieties of dimension 2 and
  is connected.
\end{theorem}

As said in the introduction, we strongly believe that the character variety of the figure 
eight knot complement does not contain any other irreducible component with irreducible 
representations. A proof should appear elsewhere.

In order to prove the theorem, we will compute 3 affine varieties $D_{1}
,D_{2} ,D_{3}$, 3 birational maps $\pi_{1} , \pi_{2} , \pi_{3}$ defined
everywhere on the $D_{i}$'s and $3$ polynomials $H_{1} ,H_{2} ,H_{3}$ in
$\mathbb{Q} [ X_{1} ,X_{2} ,X_{3} ]$ such that :
\begin{itemize}
  \item (a) $\tmop{Defor} ( M, \mathcal{T} ) = \cup_{i=1 \ldots 3} D_{i}$.
  
  \item (b) $V ( H_{i} ) = \overline{\pi_{i} ( D_{i} )} ,i=1 \ldots 3$.
  
  \item (c) $\pi_{i} ( D_{i} ) \subset V ( H_{i} )$ is smooth and $\pi_1$ realizes a homemorphism on its image.  
  \item (d) $V ( H_{i} )$ is irreducible in $\mathbb{C}^{3}$  or equivalently $H_{i}$ is
  irreducible in $\mathbb{C} [ X_{1} ,X_{2} ,X_{3} ]$.
  
  \item (e) $D_{i} \neq D_{j}$ \ for $i \neq j$ and $i,j=1 \ldots 3$.
    \item (f) $\cap_{i=1 \ldots 3} D_{i} \neq \emptyset$.
\end{itemize}
Under assumptions (a) to (f), $\pi_{i} ( D_{i} )$ is a  smooth Zariski-open subset
in the affine irreducible algebraic variety $V ( H_{i} )$, it is thus
connected -- see \cite[4.16]{Mumford}. As $\pi_{i}$ is  is a homeomorphism on its image, $D_{i}$ is also smooth and
connected. As $\cap_{i=1 \ldots 3} D_{i} \neq \emptyset$, then $\tmop{Defor}(M,\mathcal T)$
is also connected. So the theorem follows from those properties.

Let us describe how to compute $D_{i} , \pi_{i} ,H_{i} ,i=1 \ldots 3$.

\subsection{Computational tools}

The main computational tool for computing these objects is Gr\"obner basis. We recall briefly 
that a Gr\"obner basis of an ideal $\mathcal I$ in $\mathbb Q[X_1,\ldots,X_n]$ is a set of 
generators of $\mathcal I$ such that, for any polynomial $P$, there is a unique preferred 
polynomial, denoted by $P\,\; \mathrm{mod}\, \;\mathcal I$, congruent to $P$ modulo $\mathcal I$. 
This polynomial is called the \emph{normalform}. For more details, we refer to \cite{CLO1997}. 
Note that such a basis is uniquely defined once chosen an admissible ordering on the monomials. Once a Gröbner basis is known, one can compute the associated Hilbert polynomial and then the (Hilbert) dimension and the (Hilbert) degree \cite[Section 9.3]{CLO1997}.

An interesting ordering is what is called \emph{elimination ordering}. If $\mathcal I$ is seen 
as an ideal of $\mathbb Q[Y_1,\ldots, Y_k][X_1,\ldots,X_r]$,  an elimination ordering verifies 
$Y_i<X_j$ for all $i$, $j$. It allows a straightforward computation of a Gröbner basis for $$
\mathcal I'=\mathcal I\cap \mathbb Q[Y_1,\ldots,Y_k],$$ i.e. elimination of variables \cite[Ch. 
3]{CLO1997}. Indeed, the set of zeroes of $V(\mathcal I')$ is the Zariski-closure of the 
projection of $V(\mathcal I)$ on the first variables \cite[Section 4.4]{CLO1997}.

We will often eliminate some variables, especially when we may express a variable by a 
rational expression in the other ones. But, the projections might create spurious components when 
taking the Zariski closures: we work with constructible sets rather than algebraic varieties. As our computations are on the edge of what can be actually done, 
we want to avoid these components. So we will make a frequent use of saturation. Indeed, for 
an ideal $\mathcal I$ in $\mathbb Q[X_1,\ldots,X_n]$ and a  polynomial
$f$, the saturation of $\mathcal I$ by $f$ is the ideal $\mathcal I:f^\infty$,
which consist in the set of polynomials $P$ such that for some $k$,
the polynomial $Pf^k$ belongs to $\mathcal I$. Geometrically, when
$\mathcal I$, the zeroes set of $\mathcal I:f^\infty$ is the Zariski closure of $V(\mathcal I)\setminus V(f)$. The ideal $\mathcal I:f^\infty$ is computed by taking a Gröbner basis of $(\mathcal I+<tf-1>)\cap \mathbb Q[X_1,\ldots ,X_n]$, where $t$ is an independant variable -- thanks to an elimination ordering.
 
Let us formalize the operations we need and explain their scope. Here is a list of algorithms we will use later on:
\begin{itemize}
  \item [$\alpha$] Computing an ideal, in the sense of computing a Gröbner basis.
  
  \item [$\beta$] Saturating an ideal $I \in \mathbb{Q} [ X_{1} , \ldots
  ,X_{n} ]$ by  $f \in \mathbb{Q} [ X_{1} , \ldots ,X_{n}
  ]$, in the sense of computing a Gröbner basis of the ideal $I:f^{\infty}$. We
  naturally extend this algorithm for saturating an ideal $I$ by a set of
  polynomials $\{ f_{1} , \ldots ,f_{l} \}$ in order to compute $I:
  \left( \prod_{i=1 \ldots l} f_{i} \right)^{\infty}$. This is simply achieved by iteratively 
  saturating by $f_1$, then $f_2$, etc \cite[Ch. 4, prop. 10]{CLO1997}.
  
  \item [$\gamma$] Factorizing a polynomial $p \in \mathbb{Q} [ X_{1} ,
  \ldots ,X_{n} ]$ in $\mathbb{Q} [ X_{1} , \ldots ,X_{n} ]$.

  \item [$\eta$] \tmtextbf{for small systems (low degree, few number of
  variables)} Computing the prime decomposition of an ideal generated by the
  equations, removing inclusions. For those who are not familar with such
  objects, just retain that this means computing the decomposition of an
  algebraic variety in $\mathbb{Q}$- irreducible  (defined by prime
  ideals in $\mathbb{Q} [ X_{1} , \ldots ,X_{n} ]$) and non redundant
  components.
  
  \item [$\theta$] \tmtextbf{ Not usually implemented, see later on.} Testing if a polynomial $p \in \mathbb{Q} [ X_{1} ,X_{2}
  ,X_{3} ]$ is irreducible in $\mathbb{C} [ X_{1} ,X_{2} ,X_{3} ]$ or not.
\end{itemize}
 
Note that the four first algorithms are implemented in the usual computer algebra systems, such 
as Maple or Sage. Sadly, our problem seems far beyond their scope. So a great deal of our 
method is to use the shape of our equations to simplify the problem and getting to the point we can use these routines.

However the last algorithm [$\theta$] is special: it is not implemented in usual systems -- which 
work with number fields. Moreover, it is only used at the very last step of our proof of the 
theorem. So we will be more precise later on.

\subsection{An algebraic representation of $\tmop{Defor} ( M, \mathcal{T} )$}\label{ss:equations}

Given $M$ and $\mathcal{T}$ a triangulation of $M$ with $n$ tetrahedra, the
deformation variety $\tmop{Defor} (M, \mathcal{T} )$ is defined in
$\mathbb{C}^{16n}$ by a system of internal relations for each tetrahedron and
gluing equations for adjacent tetrahedra. Moreover, using the internal
relations, we may directly express the face parameters in terms of the edge
parameters. From now on, we consider $\tmop{Defor} (M, \mathcal{T} )$ as a
subset of $\mathbb{C}^{12n}$.

The system, more precisely, is given in terms of variables $Z= (z_{ij} (T))
\in \mathbb{C}^{12n}$, where $T$ denotes a tetrahedron in the triangulation
(containing $n$ tetrahedra) and $ij$ one of its oriented edges.
\begin{itemize}
  \item $L_{e}$ (resp. $L_{f}$) denotes the set of the $2n$ edge (resp. face)
  equations;
  
  \item $L_{c}$ denotes the polynomials defining the internal relations (also
  called cross-ratio relations), of the form $z_{ik} (T)  (1-z_{ij} (T)) -1$
  and $z_{ij} (T)  (1-z_{il} (T)) -1$ for each tetrahedron;
\end{itemize}
In this terms, $\tmop{Defor} (M, \mathcal{T} )$ is the algebraic variety:
\[ \tmop{Defor} (M, \mathcal{T} ) =V ( \langle L_{e} \cup L_{f} \cup L_{c}
   \rangle ) = \{{Z} \in \mathbb{C}^{12n} |P(Z)=0,P \in L_{e} \cup
   L_{f} \cup L_{c} \} . \]
Note that $\tmop{Defor} (M, \mathcal{T} ) \subset (\mathbb{C} \setminus
\{0,1\})^{12n}$, because of the internal relations $L_{c}$. This remark will be important later on.

We now display the equations for the figure eight knot complement. The triangulated structure (see section \ref{section:intro-eight}) associated
with the figure eight knot complement is made of two tetrahedra decorated by
the coordinates $(u_{ij} )_{i \neq j \in \{1, \ldots ,4\}}$ (resp. $(v_{ij}
)_{i \neq j \in \{1, \ldots ,4\}}$) with cross-ratio relations defined as the
roots of the system $L_{c_{1}} =0$ (resp. $L_{c_{2}} =0$) :
\begin{center}
  $L_{c_{1}} = \left\{ \begin{array}{c}
    u_{12} u_{13} -u_{13} +1\\
    u_{12} u_{14} -u_{12} +1\\
    u_{21} u_{23} -u_{21} +1\\
    u_{21} u_{24} -u_{24} +1\\
    u_{31} u_{34} -u_{31} +1\\
    u_{32} u_{34} -u_{34} +1\\
    u_{41} u_{43} -u_{43} +1\\
    u_{42} u_{43} -u_{42} +1
  \end{array} \right. \begin{array}{c}
    \\
    \\
    
  \end{array}$ , \ $L_{c_{2}} = \left\{ \begin{array}{c}
    v_{12} v_{13} -v_{13} +1\\
    v_{12} v_{14} -v_{12} +1\\
    v_{21} v_{23} -v_{21} +1\\
    v_{21} v_{24} -v_{24} +1\\
    v_{31} v_{34} -v_{31} +1\\
    v_{32} v_{34} -v_{34} +1\\
    v_{41} v_{43} -v_{43} +1\\
    v_{42} v_{43} -v_{42} +1
  \end{array} \right.$ ,
\end{center}
and gluing edge (resp. face) relations defined as the roots of the system
$L_{e} =0$ (resp. $L_{f} =0$) :
\begin{center}
  $L_{e} = \left\{ \begin{array}{c}
    u_{13} u_{14} u_{23} v_{21} v_{24} v_{31} -1\\
    u_{31} u_{32} u_{41} v_{12} v_{13} v_{42} -1\\
    u_{21} u_{24} u_{34} v_{23} v_{41} v_{43} -1\\
    u_{12} u_{42} u_{43} v_{14} v_{32} v_{34} -1
  \end{array} \right.$ , $L_{f} = \left\{ \begin{array}{c}
    u_{21} u_{31} u_{41} v_{21} v_{31} v_{41} -1_{\nosymbol}\\
    u_{12} u_{32} u_{42} v_{13} v_{23} v_{43} -1\\
    u_{13} u_{23} u_{43} v_{12} v_{32} v_{42} -1\\
    u_{14} u_{24} u_{34} v_{14} v_{24} v_{34} -1
  \end{array} \right.$ .
\end{center}
We set $Z= \{ u_{ij} ,v_{i,j} ,i \neq j \in \{ 1, \ldots ,4 \} \}$.

\subsection{Overview of the computations.}

Theoretically, the theorem could be proven using a prime decomposition of the ideal 
$$I_D=\langle L_{e} \cup L_{f}\cup L_{c} \rangle.$$
Indeed, as $\mathrm{Defor}(M, \mathcal T)$ is the zeroes set of this ideal, a prime decomposition would give the three components and we would be able to check their properties.
However, our system depends on $12n$ variables (24 variables in the
case of the figure eight knot complement) and is of rather high degree
(degree $6$ for some equations for the figure eight knot complement). So we are far beyond the scope of the state of the art algorithms of prime decomposition.

We will use the shape of the equations to simplify the system. The method we present is rather general for this kind of systems up to some points (number of factors
in some polynomials). It could be applied, in principle on numerous systems
with the same properties. For example, it can replace the one used in \cite{FKR14}.
Before explaining our method, let us define one more object: \emph{forbidden points}. 

Denote by $S_{D}=L_{e} \cup L_{f}\cup L_{c}$ the set of polynomials generating the ideal $I_{D}$.
Note that if $z \in V ( I_{D} )\subset \bC^{12n}$, then none of its coordinates
  belong to $\{ 1,0 \}$ \ because of the equations induced by the
  {\emph{internal relations}} $L_c$. We hence call {\tmem{forbidden
  points}} the points with any coordinate in $\{ 0,1 \}$ and we introduce 
  $F_{D}= \{ z,z-1,z \in Z \} \subset \mathbb{Q} [ Z ]$. It may seem useless at this point as this is contained in the internal relations.
  But recall that we will eliminate a lot of variables through projections. During these projections, as said before, some spurious components (i.e. all included in forbidden points) may appear. So we will always keep the equations defining the forbidden points and saturate our ideals by these equations. This will lower the number of components and ease the computations.
  
 Here are the steps of our method: 
 
  \begin{itemize}
\item {\bf Step A: elimination.} 
We iteratively detect equations in $S_{D}$ which are affine in a variable, i.e.  of the form 
$$D_{z} ( X )  z_{} -N_{z} ( X )\textrm{ with }X \subset Z \setminus \{ z \},$$ 
and with the additional condition that $D_{z} ( X )$
  cancels only at forbidden points. We can then substitute $z$ by $\frac{N_{z} ( X )}{D_{z}
  ( X )}$ in the other polynomials generating $I_{D}$, take the numerators of
  these equations and remove their factors (algorithm [$\gamma$]) that cancel
  exclusively at forbidden points (and thus never cancel on $\tmop{Defor} ( M,
  \mathcal{T} )$). This gives a new ideal with one less variable and one less equation.
  
At the end of the process, denote by $S_{Y} \subset \mathbb{Q} [ Y ]$ with $Y \subset Z$ the
  remaining equations. We keep in mind, through the set $R_{Y}
  \subset \mathbb{Q} [ Z ]$, the equations used for the substitutions. Of course, by substituting in $R_Y$, we may assume that all these equations express a substituted variable in terms of the variables in $Y$, i.e. they have the form: $D_{z} ( Y ) z-N_{z} ( Y )$, with $z
  \in Z \setminus Y$.
  
We also want to keep in mind the forbidden points. So denote by $F_{Y} \in \mathbb{Q} [ Y ]$
 the list of prime factors that appear in $F_{D}$ after having performed the same substitution.
  
We can interpret the relations $R_{Y} =0$ as a projection $\pi_{Y} : \tmop{Defor} (
  M, \mathcal{T} ) \rightarrow V ( \langle S_{Y} \rangle )$ which is not
  necessarily surjective. However, the complement of the image $V ( \langle S_{Y}
  \rangle ) \setminus \pi_{Y} ( \tmop{Defor} ( M, \mathcal{T} ) )$ only contains
  forbidden points which might have been introduced by the use of
  rational fractions in the substitutions.
  
  So we compute (algorithm $[ \beta ]$) the ideal $I_{Y}$, which is the saturation of $\langle S_{Y}  \rangle$ by $F_{Y}$ in order to get $V ( I_{Y} ) = \overline{\pi_{Y} (
  \tmop{Defor} )}$. The union of the zeroes of the polynomials $F_{Y}$ defines
  the {\tmem{forbidden points}} of $V ( I_{Y} )$ which are images by $\pi_{Y}$
  of the forbidden points of $\mathbb{C}^{12n}$.
  
Putting this in a diagram, we have :
\begin{center}
  \begin{tabular}{ccc}
    $\mathbb{C}^{12n}$ & $\rightarrow$ & $\mathbb{C}^{\sharp Y}$\\
    $\tmop{Defor} ( M, \mathcal{T} )$ & $\approx_{\pi_{Y}}$ & $V ( I_{Y} )
    \setminus ( \cup_{f \in F_{Y}} V ( f ) )$\\
    $( x,y )$ & $\rightarrow$ & $y$\\
    $\left( \frac{N_{x} ( y )}{D_{x} ( y )} ,y \right)$ & $\leftarrow$ & $y.$
  \end{tabular}
\end{center}
  \item {\bf Step B: splitting.} In $I_{Y}$, we find a polynomial
  $P_{Y} \in \mathbb{Q} [ Y ]$ that factorizes (algorithm [$\gamma$]) into
  $k$ factors in $\mathbb{Q} [ Y ]$ : $P_{Y} = \prod_{i=1 \ldots k} P_{Y,i}$.
  
  By construction, as $\tmop{Defor} ( M, \mathcal{T} )$ does not contain
  forbidden points, we then get $\tmop{Defor} ( M, \mathcal{T} ) = \cup_{i=1
  \ldots k} V ( J_{Y,i} )$, where $J_{Y,i} \subset \mathbb{Q} [ Z ]$ is the
  saturation of $I_{D} + \langle P_{Y,i} \rangle$ by $F_{D}$.
  
  We can even further saturate our ideals, in order to avoid the same component appearing twice. We  perform cascading
  saturations  for some permutation $\sigma$ of the indices: we iteratively compute
  (algorithm $[ \beta ]$)  the ideal $I_{Y, \sigma ( i )} =I_{Y} + \langle
  P_{Y, \sigma ( i )} \rangle$ and saturate this ideal by $F_{Y}$ as well as
  with $\{ P_{Y, \sigma ( 1 )} , \ldots ,P_{Y, \sigma ( i-1 )} \}$ \ in order
  to remove possible components that are already in $I_{Y, \sigma ( 1 )} ,
  \ldots ,I_{Y, \sigma ( i-1 )}$. The efficiency of the process strongly
  depends on the chosen permutation.
  
At the end of this step, one can set $D_{i} =V ( I_{Y,i} +I_{D} )$. By construction, $D_i$ is the inverse image by $\pi_i$ of the zeroes set of $I_{Y,i}$. And we get the decomposition: 
$$\tmop{Defor} ( M,
  \mathcal{T} ) = \cup_{i=1 \ldots k} D_{i}.$$
  
  \item {\bf Step C: further elimination.} The ideals $I_{Y,i}$ are still too big to perform 
  directly a prime decomposition. So, for each $I_{Y,i}$, we apply once again step A. We may 
  even try to iteratively apply Steps A and B and so on, but in our case, a single additional 
  pass of Step A turns out to be enough. So we explain only this single pass.
  
  Indeed, for each $I_{Y,i}$, we find a polynomial that is affine in
  one of the variables and, moreover, can be written $a_{i} ( z ) y+b_{i} ( z
  )$ with $z \in Y_{i} =Y \setminus \{ y \}$ and $a_{i} ( z^{} )$ only vanishing at forbidden points. We stress that 
  the suitable variable $y$ to eliminate depends on the different ideals $I_{Y,i}$.
  
  The result of Step A is then an ideal $I_{Y_{i}} \subset \mathbb{Q}
  [ Y_{i} ]$ together with a projection $\pi_{Y_{i}} :V ( I_{Y,i} ) \rightarrow V (
  I_{Y_{i}} )$ and a set of polynomials
  $F_{Y_{i}}$ (defining the {\tmem{forbidden points}} of $V ( I_{Y_{i}} )$)
  such that we get: $$\overline{\pi_{Y_{i}} ( V ( I_{Y,i} ) )} =V ( I_{Y_{i}} ).$$
  
  Setting $\pi_{i} = \pi_{Y_{i}} \circ \pi_{Y}$, we may write
  $\overline{\pi_{i} ( D_{i} )} =V ( I_{Y_{i}} )$ and get the diagram similar to the previous 
  one:
 
\begin{center}
    \begin{tabular}{ccc}
      $\mathbb{C}^{12n}$ & $\rightarrow$ & $\mathbb{C}^{\sharp Y-1}$ \\

      $D_{i}$ & $\approx_{\pi_{Y_i}}$ & $V ( I_{Y_i} ) \setminus ( \cup_{f \in
      F_{Y_i}} V ( f ) )$ \\
%
%
    \end{tabular}
\end{center}
%
%
%
Beware that the different projections $\pi_{Y_i}$ are only defined on their respective component $D_i$ and do not land in the same space, as the last eliminated variable may differ for different components.
  
  \item {\bf Step D: Prime decomposition.} At this stage, each ideal $I_{Y_{i}}$ belongs to 
  the scope of a Prime Decomposition algorithm (named here algorithm [$\eta$]). Doing so, we 
  extract the $\mathbb{Q}$-irreducible components or just a check that the $I_{Y_{i}}$ 
  is prime. We denote by $I_{Y_i,k}$ the different prime ideals appearing in the decomposition of $I_{Y_i}$. 
\end{itemize}

The Steps $A$ to $D$ are general functions that also have been applied
for other computations, for example for computing unipotent solutions for many
cases of varieties with a triangulation involving up to $4$ tetrahedra (\cite{FKR14}).

At this point, we have a decomposition of $\mathrm{Defor}(M,\mathcal T)$ in components, for which we have a prime decomposition. Proving the theorem for this decomposition is now just a matter of checking properties and here the method is very specific to the figure eight knot complement. Other manifolds $M$ may give other behaviors. In particular, it turns out in our case that, at this step, each ideal $I_{Y_i}$ has a component principal and prime and maybe another component of lesser dimension and which is redundant:
\begin{itemize}
  \item {\bf Specific Step $1$: irreducibility and smoothness} Check that for each $i$ there is a unique prime ideal of maximal dimension among the prime factors. We denote here by $I_{Y_i,1}$ this prime factor. Check that $\mathrm{Defor}(M,\mathcal T)$ is included in the union of the $\pi_i^{-1}(I_{Y_i,1})$. In other terms, the components of lesser dimension are redundant. Then check that each $I_{Y_{i} ,1} \subset \mathbb{Q} [ Y_{i} ]$ is  a  principal ideal of $\mathbb{Q} [ Y_{i} ]$ (the output of algorithm
  [$\alpha$] is a unique polynomial), and, if so, denote by $H_{i} \in
  \mathbb{Q} [ Y_{i} ]$ its generator.
  
  Eventually verify (algorithm [$\gamma$]) that $H_{i}$ is irreducible in
  $\mathbb{C} [ Y_{i} ]$ and check that the singular points of each $H_{i}$
  are exclusively forbidden points.
  
  \item {\bf Specific Step $2$: connectedness.} Let $I_{i}$ be the ideal generated by $I_{Y_i,1}$ and $I_D$. It verifies that $V ( I_{i} ) =D_{i}$. We compute (algorithm [$\alpha$]) the ideal $I_{1} + \ldots +I_{k}$. And we check if is not trivial.
\end{itemize}
If the results of these two specific steps are indeed what is announced, then the theorem is easily proven: $\pi_{i} ( D_{i} )$ is a dense open
  set in $V ( H_{i} )$ which is a smooth affine (irreducible) variety so that
  $\pi_{i} ( D_{i} )$ and thus $D_{i}$ are smooth and connected. Moreover,  as the ideal $I_{1} + \ldots +I_{k}$ is not trivial, we get that $D_{1} \cap \ldots \cap D_{k}$ is a non empty 
  set of points. As each $D_i$ is connected, we conclude that $\tmop{Defor} ( M, \mathcal{T} ) = \cup_{i=1 \ldots k} D_{i}$ is connected.

\subsection{Explicit computations for the Figure Eight Knot complement.}\label{results}

We now review the different steps in the case of the figure eight knot complement. We indicate the choices that can be done to ensure the completion of our method. Beware that the feasibility of the computations depends on these choices. For example, there are choices of the eliminated variables in Step A that lead to unfeasible computations. Let us mention that the reader may find the Maple worksheet where this computation is implemented: see \cite{Worksheet}.

\subsubsection{{\bf Step A}}

We apply step A on $S_{D} =L_{c_{1}} \cup L_{c_{2}} \cup L_{e} \cup
L_{f}$. With the set $Y= \{ v_{1,4} ,u_{4,3} ,v_{4,3} ,v_{3,4} \}$, one get the following expressions for $\pi_Y$:
{\small \begin{center}
  $\pi_{Y} : \left\{ \begin{array}{l}
    u_{12} =1- \frac{1}{u_{13}} ,u_{13} =1- \frac{1}{u_{14}} ,u_{14} =1-
    \frac{(1-v_{34} )(1-u_{24} )}{1-v_{24}} ,u_{21} = \frac{1}{1-u_{23}} ,\\
    u_{23} = \frac{1}{1-u_{24}} ,u_{24} =1+v_{14} (1-v_{24} )
    \frac{u_{43}}{1-u_{43}} ,u_{31} = \frac{1}{1-u_{34}} ,u_{32} =1-
    \frac{1}{u_{34}} ,\\
    u_{34} =- \frac{1-u_{43}}{u_{43} v_{14} (1-v_{43} )} ,u_{41} =1-
    \frac{1}{u_{43}} ,u_{42} = \frac{1}{1-u_{43}} ,\\
    v_{12} = \frac{1}{1-v_{14}} ,v_{13} =1- \frac{1}{v_{14}} ,v_{21} =1-
    \frac{1}{v_{24}} ,v_{23} = \frac{1}{1-v_{24}} ,\\
    v_{24} =1- \frac{v_{43} (1-v_{14} )}{u_{43} v_{14}^{2} (1-v_{34} )} -
    \frac{v_{43} (1-v_{43} )(1-v_{14} )}{(1-u_{43} )v_{14} (1-v_{34} )} ,\\
    v_{31} = \frac{1}{1-v_{34}} ,v_{32} =1- \frac{1}{v_{34}} ,v_{41} =1-
    \frac{1}{v_{43}} ,v_{42} = \frac{1}{1-v_{43}} .
  \end{array} \right.$
\end{center}}
Moreover the ideal $I_{Y}$ is generated by $3$ polynomials in $\mathbb{Q} [
Y ]$ of maximal degree $5$ in each variable, which are too large to be printed
here. Neither do we print $F_{Y}$ which is quite large and easy to compute using a
simple substitution and some factorizations.

Note that $\pi_{Y}$ is not simplified here, the expression having only
variables of $Y$ in the left sizes of the equations being a little too large
to be printed in this article.

\subsubsection{{\bf Step B}}

We are looking for a polynomial in $I_Y$ that may be factorized. As said before, $I_Y$ is generated by $3$ polynomials. Computing the gcd in $\mathbb{Q} [ v_{1,4} ,u_{4,3} ,v_{4,3} ]$ of all the resultants of 2 of these $3$ polynomials wrt $v_{3,4}$ we then obtain a
polynomial that belongs to $I_{Y}$ and factorizes into $3$ factors. We use this polynomial as $P_{Y}$ and get the three factors:
\begin{itemize}
  \item {\small{$P_{Y,1} =-u_{4,3} v_{1,4}^{2} v_{4,3}^{2} +u_{4,3}^{2} v_{1,4}
  v_{4,3} +u_{4,3} v_{1,4}^{2} v_{4,3} +u_{4,3} v_{1,4} v_{4,3}^{2}
  -u_{4,3}^{2} v_{1,4} -3 \hspace{0.17em} u_{4,3} v_{1,4} v_{4,3} +v_{1,4}
  u_{4,3} +u_{4,3} v_{4,3} +v_{1,4} v_{4,3} -v_{4,3}$}}
  
  \item {\small{$P_{Y,2} =-v_{1,4}^{2} v_{4,3}^{2} +2 \hspace{0.17em}
  v_{1,4}^{2} v_{4,3} +v_{1,4} v_{4,3}^{2} +u_{4,3} v_{4,3} -v_{1,4}^{2}
  -v_{1,4} v_{4,3} -v_{4,3}$}}
  
  \item {\small{$P_{Y,3} =u_{4,3}^{5} v_{1,4}^{4} v_{4,3}^{2} +u_{4,3}^{4}
  v_{1,4}^{4} v_{4,3}^{3} +u_{4,3}^{3} v_{1,4}^{4} v_{4,3}^{4} -u_{4,3}^{5}
  v_{1,4}^{4} v_{4,3} -3 \hspace{0.17em} u_{4,3}^{4} v_{1,4}^{4} v_{4,3}^{2}
  -u_{4,3}^{4} v_{1,4}^{3} v_{4,3}^{3} -3 \hspace{0.17em} u_{4,3}^{3}
  v_{1,4}^{4} v_{4,3}^{3} -2 \hspace{0.17em} v_{1,4}^{3} u_{4,3}^{3}
  v_{4,3}^{4} +2 \hspace{0.17em} v_{1,4}^{4} u_{4,3}^{4} v_{4,3} +2
  \hspace{0.17em} v_{1,4}^{3} u_{4,3}^{4} v_{4,3}^{2} +3 \hspace{0.17em}
  u_{4,3}^{3} v_{4,3}^{2} v_{1,4}^{4} +6 \hspace{0.17em} v_{1,4}^{3}
  u_{4,3}^{3} v_{4,3}^{3} +u_{4,3}^{3} v_{1,4}^{2} v_{4,3}^{4} -u_{4,3}^{4}
  v_{1,4}^{3} v_{4,3} -3 \hspace{0.17em} u_{4,3}^{4} v_{1,4}^{2} v_{4,3}^{2}
  -u_{4,3}^{3} v_{1,4}^{4} v_{4,3} -6 \hspace{0.17em} u_{4,3}^{3} v_{1,4}^{3}
  v_{4,3}^{2} -5 \hspace{0.17em} u_{4,3}^{3} v_{1,4}^{2} v_{4,3}^{3}
  -u_{4,3}^{2} v_{1,4}^{3} v_{4,3}^{3} +3 \hspace{0.17em} v_{1,4}^{2}
  u_{4,3}^{4} v_{4,3} +2 \hspace{0.17em} v_{1,4}^{3} u_{4,3}^{3} v_{4,3} +10
  \hspace{0.17em} u_{4,3}^{3} v_{4,3}^{2} v_{1,4}^{2} +2 \hspace{0.17em}
  u_{4,3}^{3} v_{1,4} v_{4,3}^{3} +2 \hspace{0.17em} u_{4,3}^{2} v_{1,4}^{3}
  v_{4,3}^{2} +3 \hspace{0.17em} u_{4,3}^{2} v_{1,4}^{2} v_{4,3}^{3}
  -u_{4,3}^{4} v_{1,4} v_{4,3} -6 \hspace{0.17em} v_{1,4}^{2} u_{4,3}^{3}
  v_{4,3} -2 \hspace{0.17em} u_{4,3}^{3} v_{1,4} v_{4,3}^{2} -v_{1,4}^{3}
  u_{4,3}^{2} v_{4,3} -6 \hspace{0.17em} v_{1,4}^{2} u_{4,3}^{2} v_{4,3}^{2}
  -2 \hspace{0.17em} u_{4,3}^{2} v_{1,4} v_{4,3}^{3} +2 \hspace{0.17em}
  u_{4,3}^{3} v_{1,4} v_{4,3} +u_{4,3}^{3} v_{4,3}^{2} +3 \hspace{0.17em}
  u_{4,3}^{2} v_{1,4}^{2} v_{4,3} +2 \hspace{0.17em} u_{4,3}^{2} v_{1,4}
  v_{4,3}^{2} +u_{4,3}^{4} +u_{4,3}^{3} v_{4,3} -u_{4,3}^{2} v_{1,4} v_{4,3}
  -2 \hspace{0.17em} u_{4,3}^{2} v_{4,3}^{2} -4 \hspace{0.17em} u_{4,3}^{3} -3
  \hspace{0.17em} u_{4,3}^{2} v_{4,3} +u_{4,3} v_{4,3}^{2} +6 \hspace{0.17em}
  u_{4,3}^{2} +3 \hspace{0.17em} u_{4,3} v_{4,3} -4 \hspace{0.17em} u_{4,3}
  -v_{4,3} +1$}}
\end{itemize}
As mentioned in the above section, the efficiency of the ``cascading
saturations'' of step [B] strongly depends on the numbering of these factors.
In the present case, a favorable permutation of the indices is $\sigma ( 1 ) =2,
\sigma ( 2 ) =1$ and $\sigma ( 3 ) =3$, other permutations might drive to
infeasible computations.

We thus set $I_{Y,2} =I_{Y} + \langle P_{Y,2} \rangle$ and saturate this ideal
by $F_{Y}$, we then compute $I_{Y,1} =I_{Y} + \langle P_{Y,1} \rangle$ and
saturate this ideal by $F_{Y} \cup \{ P_{Y,2} \}$ and finally compute $I_{Y,3}
=I_{Y} + \langle P_{Y,3} \rangle$ and saturate this ideal by $F_{Y} \cup \{
P_{Y,2} ,P_{Y,3} \}$.

\subsubsection{{\bf Step C}}

For each of these three ideals, we want to eliminate one more variable. Here are the affine polynomials we use (checking each time that the leading coefficient only vanishes at forbidden points):
\begin{enumerate}
\item In $I_{Y,1}$, we find the polynomial $v_{3,4} -v_{4,3}$ and eliminate $v_{3,4}$.
\item
In $I_{Y,2}$, the polynomial $P_{Y,2}$ displayed above is affine in $u_{4,3}$ with leading coefficient $v_{4,3}$. We eliminate $u_{4,3}$.
\item In $I_{Y,3}$, we find the following polynomial: 
$$\left(v_{1,4}^{2}v_{3,4}(v_{3,4}-1)\right) u_{4,3} -v_{1,4} v_{3,4} +1.$$ We eliminate $u_{4,3}$.
\end{enumerate}
We then define the three projections $\pi_{1}$, $\pi_{2}$ and  $\pi_{3}$ and compute (algorithm [$\alpha$]) the ideals
$I_{Y_{1}} ,I_{Y_{2}}$ and $I_{Y_{3}}$ by
substituting in $I_{Y}$ and saturating by respectively $F_{Y} ,F_{Y} \cup \{
P_{Y,1} \} ,F_{Y} \cup \{ P_{Y,1} ,P_{Y,2} \}$.

\subsubsection{{\bf Step D}}

We may now compute the prime decompositions of the three ideals. The output is as follows: the above ideal
$I_{Y_{1}}$, $I_{Y_{2}}$ are not prime and $I_{Y_{3}}$ is prime. After having
removed the embedded components, it appears that $I_{Y} +I_{Y_{1}} =I_{Y_{1}
,1} \cap I_{Y_{1} ,2}$ and $I_{Y} +I_{Y_{2}} =I_{Y_{2} ,1} \cap I_{Y_{2} ,2}$. In this decomposition,
$I_{Y_{1} ,1}$ and $I_{Y_{2} ,1}$ are prime ideals of dimension $1$, whereas $I_{Y_{1} ,2}$, $I_{Y_{2} ,2}$ and $I_{Y} +I_{Y_{3}}$ are prime ideals of dimension $2$.

\subsubsection{{\bf Specific Step $1$}}
We check that $I_{Y_{1} ,1} \cap ( I_{Y} +I_{Y_{3}} ) =I_{Y_{1} ,1}$ and
$I_{Y_{2} ,1} \cap ( I_{Y} +I_{Y_{3}} )$ so that variety
$$V ( I_{Y} ) =V ( I_{Y} +I_{Y_{1} ,2} ) \cup V ( I_{Y} +I_{Y_{2} ,2}
  ) \cup V ( I_{Y} +I_{Y_{3}} )$$ is the union of $3$ $\mathbb{Q}$-irreducible
  algebraic varieties of dimension $2$. 
  
 We check that the three ideals $I_{Y_1,2}$, $I_{Y_2,2}$ and $I_{Y_3}$ are
  principal ideals each defined by a unique polyomial $H_{i}
\in \mathbb{Q} [ Y_{i} ]$. And we additionally check that the $H_i$'s are irreducible in $\mathbb{C} [ Y_{i} ]$ (algorithm [$\theta$]). Here are the polynomials $H_i$'s:
\begin{itemize}
  \item {\small{$H_{1} =u_{4,3} v_{1,4}^{2} v_{4,3}^{2} -u_{4,3}^{2} v_{1,4}
  v_{4,3} -u_{4,3} v_{1,4}^{2} v_{4,3} -u_{4,3} v_{1,4} v_{4,3}^{2}
  +u_{4,3}^{2} v_{1,4} +3 \hspace{0.17em} u_{4,3} v_{1,4} v_{4,3} -v_{1,4}
  u_{4,3} -u_{4,3} v_{4,3} -v_{1,4} v_{4,3} +v_{4,3}$}}
  
  \item {\small{$H_{2} =v_{1,4}^{5} v_{3,4}^{2} v_{4,3}^{4} -4 \hspace{0.17em}
  v_{4,3}^{3} v_{1,4}^{5} v_{3,4}^{2} -2 \hspace{0.17em} v_{1,4}^{5} v_{3,4}
  v_{4,3}^{4} -3 \hspace{0.17em} v_{4,3}^{4} v_{1,4}^{4} v_{3,4}^{2} +6
  \hspace{0.17em} v_{4,3}^{2} v_{1,4}^{5} v_{3,4}^{2} +8 \hspace{0.17em}
  v_{1,4}^{5} v_{3,4} v_{4,3}^{3} +v_{1,4}^{5} v_{4,3}^{4} +9 \hspace{0.17em}
  v_{4,3}^{3} v_{1,4}^{4} v_{3,4}^{2} +5 \hspace{0.17em} v_{1,4}^{4} v_{3,4}
  v_{4,3}^{4} +3 \hspace{0.17em} v_{1,4}^{3} v_{3,4}^{2} v_{4,3}^{4}
  -v_{1,4}^{3} v_{3,4} v_{4,3}^{5} -4 \hspace{0.17em} v_{1,4}^{5} v_{3,4}^{2}
  v_{4,3} -12 \hspace{0.17em} v_{1,4}^{5} v_{3,4} v_{4,3}^{2} -4
  \hspace{0.17em} v_{1,4}^{5} v_{4,3}^{3} -9 \hspace{0.17em} v_{3,4}^{2}
  v_{1,4}^{4} v_{4,3}^{2} -15 \hspace{0.17em} v_{3,4} v_{4,3}^{3} v_{1,4}^{4}
  -2 \hspace{0.17em} v_{1,4}^{4} v_{4,3}^{4} -5 \hspace{0.17em} v_{1,4}^{3}
  v_{3,4}^{2} v_{4,3}^{3} -v_{1,4}^{3} v_{4,3}^{4} v_{3,4} -v_{1,4}^{2}
  v_{3,4}^{2} v_{4,3}^{4} +3 \hspace{0.17em} v_{1,4}^{2} v_{4,3}^{5} v_{3,4}
  +v_{3,4}^{2} v_{1,4}^{5} +8 \hspace{0.17em} v_{1,4}^{5} v_{3,4} v_{4,3} +6
  \hspace{0.17em} v_{1,4}^{5} v_{4,3}^{2} +3 \hspace{0.17em} v_{3,4}^{2}
  v_{1,4}^{4} v_{4,3} +15 \hspace{0.17em} v_{1,4}^{4} v_{3,4} v_{4,3}^{2} +6
  \hspace{0.17em} v_{1,4}^{4} v_{4,3}^{3} +v_{1,4}^{3} v_{3,4}^{2} v_{4,3}^{2}
  +2 \hspace{0.17em} v_{4,3}^{3} v_{1,4}^{3} v_{3,4} +v_{1,4}^{3} v_{4,3}^{4}
  -v_{1,4}^{2} v_{3,4}^{2} v_{4,3}^{3} -7 \hspace{0.17em} v_{4,3}^{4}
  v_{1,4}^{2} v_{3,4} -3 \hspace{0.17em} v_{1,4} v_{3,4} v_{4,3}^{5} -2
  \hspace{0.17em} v_{1,4}^{5} v_{3,4} -4 \hspace{0.17em} v_{1,4}^{5} v_{4,3}
  -5 \hspace{0.17em} v_{4,3} v_{1,4}^{4} v_{3,4} -6 \hspace{0.17em}
  v_{1,4}^{4} v_{4,3}^{2} +v_{1,4}^{3} v_{3,4}^{2} v_{4,3} +3 \hspace{0.17em}
  v_{3,4} v_{1,4}^{3} v_{4,3}^{2} +2 \hspace{0.17em} v_{1,4}^{2} v_{3,4}^{2}
  v_{4,3}^{2} +10 \hspace{0.17em} v_{1,4}^{2} v_{3,4} v_{4,3}^{3} +v_{1,4}
  v_{3,4}^{2} v_{4,3}^{3} +7 \hspace{0.17em} v_{1,4} v_{3,4} v_{4,3}^{4}
  +v_{3,4} v_{4,3}^{5} +v_{1,4}^{5} +2 \hspace{0.17em} v_{1,4}^{4} v_{4,3} -3
  \hspace{0.17em} v_{1,4}^{3} v_{3,4} v_{4,3} -3 \hspace{0.17em} v_{1,4}^{3}
  v_{4,3}^{2} -6 \hspace{0.17em} v_{1,4}^{2} v_{3,4} v_{4,3}^{2} -2
  \hspace{0.17em} v_{1,4}^{2} v_{4,3}^{3} -6 \hspace{0.17em} v_{1,4} v_{3,4}
  v_{4,3}^{3} -2 \hspace{0.17em} v_{3,4} v_{4,3}^{4} +2 \hspace{0.17em}
  v_{1,4}^{3} v_{4,3} +2 \hspace{0.17em} v_{1,4}^{2} v_{4,3}^{2} +v_{3,4}
  v_{4,3}^{3} +v_{1,4} v_{4,3}^{2}$}}
  
  \item {\small{$H_{3} =v_{1,4}^{5} v_{3,4}^{3} v_{4,3}^{2} -v_{1,4}^{5}
  v_{3,4}^{3} v_{4,3} +v_{1,4}^{4} v_{3,4}^{5} +v_{3,4}^{4} v_{4,3}
  v_{1,4}^{4} -2 \hspace{0.17em} v_{1,4}^{4} v_{3,4}^{3} v_{4,3}^{2} -3
  \hspace{0.17em} v_{3,4}^{4} v_{1,4}^{4} -v_{3,4}^{2} v_{1,4}^{4} v_{4,3}^{2}
  -v_{3,4}^{4} v_{4,3} v_{1,4}^{3} +v_{1,4}^{3} v_{3,4}^{3} v_{4,3}^{2} +3
  \hspace{0.17em} v_{1,4}^{4} v_{3,4}^{3} +2 \hspace{0.17em} v_{3,4}^{2}
  v_{1,4}^{4} v_{4,3} -2 \hspace{0.17em} v_{1,4}^{3} v_{3,4}^{4} +2
  \hspace{0.17em} v_{1,4}^{3} v_{3,4}^{2} v_{4,3}^{2} -v_{1,4}^{4} v_{3,4}^{2}
  +4 \hspace{0.17em} v_{1,4}^{3} v_{3,4}^{3} -2 \hspace{0.17em} v_{1,4}^{3}
  v_{3,4}^{2} v_{4,3} +v_{1,4}^{2} v_{3,4}^{3} v_{4,3} -v_{1,4}^{2}
  v_{3,4}^{2} v_{4,3}^{2} -2 \hspace{0.17em} v_{1,4}^{3} v_{3,4}^{2} +3
  \hspace{0.17em} v_{1,4}^{2} v_{3,4}^{3} +v_{1,4}^{2} v_{3,4}^{2} v_{4,3} -5
  \hspace{0.17em} v_{1,4}^{2} v_{3,4}^{2} -v_{1,4}^{2} v_{3,4} v_{4,3}
  -v_{1,4} v_{3,4}^{2} v_{4,3} +2 \hspace{0.17em} v_{1,4}^{2} v_{3,4} -2
  \hspace{0.17em} v_{1,4} v_{3,4}^{2} +v_{1,4} v_{3,4} v_{4,3} +2
  \hspace{0.17em} v_{1,4} v_{3,4} +v_{3,4} -1$.}}
\end{itemize}
Showing that, in addition, the $H_{i}$ are not singular outside the forbidden points is a 
straightforward application of the algorithms $[ \alpha ]$ and $[\gamma]$: the ideal generated 
by the partial derivatives of $H_i$, once saturated by $F_{Y_i}$, is trivial.

\subsubsection{{\bf Specific Step $2$}}

We complete the computations by a direct application of algorithm $[
\alpha ]$ in order to show that $I_{Y_{1}} +I_{Y_{2}}
+I_{Y_{3}}$ is non trivial and zero-dimensional.

\subsection{Details on algorithm [$\theta$]}

In this section we give more details on the algorithm [$\theta$] which is not available
in usual computer algebra systems.

Algorithm $[ \theta ]$ (factorization in $\mathbb{C} [ X_{1} , \ldots ,X_{n}
]$) is not present in most of computer algebra systems since they
don't know, in general, how to perform exact operations in
$\mathbb{C}$. On the other hand, most computer algebra systems implements algorithms for factorizing in
$\mathbb{Q} [ X, \ldots ,X_{n} ]$ or even $\mathbb{Q} ( \alpha ) [ X_{1} ,
\ldots ,X_{n} ]$ where $\alpha$ is any algebraic number defined by a
$\mathbb{Q}$-irreducible polynomial.

A first way to avoid this problem is to use the specific form of the $H_i$'s. Indeed, all three of them are quadratic in one variable. An elementary and careful examination, which mostly reduces to computing discriminants and checking they are not squares, leads to the proof those polynomials are irreducible over $\bC$.

More generally, a method for factorizing a polynomial $P$ over $\mathbb C$ is to find a number field
$\mathbb{Q} ( \alpha )$ such that all its factors belong to $\mathbb{Q} (
\alpha ) [ X_{1} , \ldots ,X_{n} ]$. Given a suitable $\alpha$ defined by its
minimal polynomial (with coefficients in $\mathbb{Q}$) the factorization in
$\mathbb{Q} ( \alpha ) [ X_{1} , \ldots ,X_{n} ]$ can be performed by most
computer algebra systems. Some methods for finding such a suitable $\alpha$
are reviewed in \cite{Rag97}, including the algorithm implemented in our worksheet.

\section{Eigenvalues of peripheral representations}\label{s:A-var}

Given the triangulated structure associated with the figure eight knot
complement, the meridian and the longitude are computed from the Snappea
triangulation. The corresponding diagonal entries of the holonomy matrices are \cite{FKR14}:
\[ g_{m} = \left[ \begin{array}{ccc}
     \frac{u_{12} u_{21}}{u_{34} w_{24}} & \ast & \ast \\
     0 & 1 & \ast \\
     0 & 0 & \frac{v_{42}}{u_{34}}
   \end{array} \right] \]
\begin{center}
  $g_{l} = \left[ \begin{array}{ccc}
    \frac{u_{13} u_{31} v_{14} v_{23}}{v_{31} u_{42} v_{42} u_{24}} & \ast &
    \ast\\
    0 & 1 & \ast\\
    0 & 0 & \frac{u_{31} v_{13} u_{13} v_{24}}{u_{42} v_{32} u_{24} v_{41}}
  \end{array} \right]$
\end{center}
We therefore deduce the eigenvalues of the holonomies of the longitude and
meridian:
\[ L= \frac{u_{12} u_{21}}{u_{34} v_{24}} ,L^{\ast} = \frac{v_{42}}{u_{34}}
   ,M= \frac{u_{13} u_{31} v_{14} v_{23}}{v_{31} u_{42} v_{42} u_{24}}
   ,M^{\ast} = \frac{u_{31} v_{13} u_{13} v_{24}}{u_{42} v_{32} u_{24} v_{41}}
   . \]

\subsection{$A$-variety}

The study of these eigenvalues on each component of $\mathrm{Defor}(M,\mathcal T)$ can then be done by adding the following set $L_{\{l,m\}}$ of equations to the generators of, respectively, $I_{Y_{1}} ,I_{Y_{2}}$ and $I_{Y_{3}}$:
\begin{eqnarray*}
L_{\{ l,m \}} &= &\left\{ u_{34} v_{24} L-u_{12} u_{21}  , \quad v_{31} u_{42} v_{42} u_{24} M-u_{13} u_{31} v_{14} v_{23}  ,\right. \\
&&\left.  u_{34} L^{\ast} -v_{42} 
, \quad u_{42} v_{32}
u_{24} v_{41} M^{\ast} -u_{31} v_{13} u_{13} v_{24} \right\}.
\end{eqnarray*}

Eliminating all the variables but $\{ L,L^{\star} ,M,M^{\star} \}$ by means of
Gr{\"o}bner bases and similar factorization tricks as for the computation of
$\tmop{Defor} ( M, \mathcal{T} )$, we get that :
\begin{itemize}
  \item $( I_{Y_{2}} + \langle L_{l,m} \rangle ) \cap \mathbb{Q} [
  L,L^{\star} ,M,M^{\star} ] = \langle L^{3} -M,L^{\star 3} -M^{\star}
  \rangle$ or, equivalently, that over all the points of $D_{2}$, $L^{3} =M$
  and $L^{\star 3} =M^{\star}$
  
  \item $( I_{Y_{3}} + \langle L_{l,m} \rangle ) \cap \mathbb{Q} [
  L,L^{\star} ,M,M^{\star} ] = \langle L^{3} M-1,L^{\star 3} M^{\star} -1
  \rangle$ or, equivalently, that over all the points of $D_{3}$, $L^{3} =
  \frac{1}{M}$ and $L^{\star 3} = \frac{1}{M^{\star}}$
\end{itemize}
The eigenvalues over the component $D_{1}$ are more complicated to describe.
The Gr{\"o}bner basis of $( I_{Y_{1}} + \langle L_{l,m} \rangle ) \cap
\mathbb{Q} [ L,L^{\star} ,M,M^{\star} ]$, for the Degree Reverse
Lexicographic ordering is huge (made of 141 polynomials with up to 3462
terms). This might be due to the ordering itself or to the fact that the ideal
might have many complicated embedded components which can be viewed as an
artifact of the projection onto the coordinates $[L,L^{\star} ,M,M^{\star} ]$.

Using an elimination ordering with $M>L,L^{\star} ,M^{\star}$ one gets a
simpler description showing in particular that generically (say outside a
subvariety of dimension $1$), the points of that component of the $A$-variety
are in one-to-four correspondence with those of the hypersurface $H_{1}$ of
$\mathbb{C}^{3}$.

\subsection{Where are some already known representations ?}

\label{section:example-unipotent}Some representations are already known, more
precisely the boundary-unipotent ones, i.e. those for which the images of $l$
and $m$ are unipotent. These representations were found in {\cite{FKR14}}.
Algebraically, we add the conditions $L=L^{\ast} =M=M^{\ast} =1$. Once a point in $\mathrm{Defor}(M,\mathcal T)$ is known, it is an easy task to decide whether it belongs to $D_1$, $D_2$ or $D_3$: just plug the coordinates in the polynomials defining these three components and check which one vanishes.

Here are the representations found in \cite{FKR14} and where they live:
\begin{itemize}
  \item The monodromy of the complete hyperbolic structure on the complement of the figure-eight knot and its complex conjugate representation correspond to points of $D_1$. Indeed, if
  $\omega^{\pm} = \frac{1 \pm i \sqrt{3}}{2}$ is one root of $f_{1}$, then
  \[ u_{12} =u_{21} =u_{34} =u_{43} =v_{12} =v_{21} =v_{34} =v_{43} =
     \omega^{\pm} \]
  is a point of $\mathrm{Defor}(M,\mathcal T)$, corresponding to this monodromy or its complex conjugate. In other terms, $D_1$ is the geometric component of $\mathrm{Defor}(M,\mathcal T)$.
  
  \item Setting:
  \[ u_{1,2} =u_{3,4} =v_{3,4} =v_{4,3} = \bar{u}_{2,1} = \bar{u}_{4,3} =
     \bar{v}_{1,2} = \bar{v}_{2,1} = \omega^{\pm} , \]
     one gets two other points in $\mathrm{Defor}(M,\mathcal T)$.
  They correspond to a discrete representation of the fundamental
  group of the complement of knot in $\tmop{PU} (2,1)$ with faithful boundary
  holonomy. Moreover, its action on complex hyperbolic space has limit set the
  full boundary sphere {\cite{falbeleight}}. Both these point also belong to the geometric component $D_1$.

  \item Let $\gamma^{\pm} =-
  \frac{1}{2} \pm i \frac{1}{2}  \sqrt{7}$. Define
  \begin{center}
    $\begin{array}{c}
      u_{2,1} =v_{2,1} = \bar{u}_{4,3} = \bar{v}_{1,2} = \frac{5\pm i
      \sqrt{7}}{4} ,\\
      u_{1,2} =v_{2,1} = \frac{3\pm i \sqrt{7}}{8} ,v_{4,3} = \bar{v}_{3,4} =-
      \frac{1\pm i \sqrt{7}}{2}
    \end{array}$
  \end{center}
  These two points belong to $\mathrm{Defor}(M,\mathcal T)$ and correspond to representation in $\tmop{PU} (2,1)$. It turns out that they belong to $D_3$.
  \item Consider now the inverse of $\gamma^{\pm}$: $1/ \gamma^{\pm} =-
  \frac{1}{4} \pm i \frac{1}{4}  \sqrt{7}$,nd define:
  \begin{center}
    $\begin{array}{c}
      u_{1,2} = \bar{u}_{3,4} = \frac{3\pm i \sqrt{7}}{2} ,v_{4,3} =
      \bar{v}_{3,4} = \frac{5\pm i \sqrt{7}}{8} ,\\
      u_{2,1} =v_{2,1} = \bar{u}_{4,3} = \bar{v}_{1,2} = \frac{-1\pm i
      \sqrt{7}}{4}
    \end{array}$
  \end{center}  
  These two points belong to $\mathrm{Defor}(M,\mathcal T)$ and correspond to representation in $\tmop{PU} (2,1)$. It turns out that they belong to $D_2$.

\end{itemize}
All solutions were already obtained in {\cite{falbeleight}}. Moreover the solutions that belong to $D_2$ or $D_3$ correspond to spherical CR structures with
unipotent boundary holonomy of rank one {\cite{DF}}.

Moreover, in {\cite{DF}}, it was shown that there exists a non-inner
automorphism $\tau$ of the group $\Gamma_{8}$ that sends a representation
appearing in the third point above to a representation appearing in the fourth point.
It is known that, on the hyperbolic structure, this automorphism
correspond to a change of orientation, that is it acts as complex conjugation. It proves that the action of $\tau$
on the representation variety exchange the last two components and preserves
the first one. This fact is presented in a more concrete way in the following section
\ref{deformation}.

\subsubsection{Hyperbolic solutions - $A$-polynomial}

We may now look after the representations coming from $\mathrm{PGL}
(2,\mathbb{C})$, which we call hyperbolic solutions. Indeed, the hyperbolic
conditions : $u_{12} =u_{21} =u_{34} =u_{43}$, $v_{12} =v_{21} =v_{34}
=v_{43}$ on the parameters is equivalent to the fact that the holonomy
representation associated to these parameters has its image lying in
$\mathrm{Ad} ( \mathrm{PGL} (2,\mathbb{C}))$. These conditions simplify a lot
the equations and give the 1-dimensional complex variety
\[ z= \frac{1}{1-y} , \hspace{1em} x^{2} y^{2} = (1-x)  (1-y) . \]
We then obtain
\[ \begin{array}{l}
     u_{12} =u_{34} =u_{21} =u_{43} =x,\\
     u_{23} =u_{41} =u_{14} =u_{32} =1- \frac{1}{x} ,\\
     u_{13} =u_{31} =u_{24} =u_{42} = \frac{1}{1-x} ,\\
     v_{12} =v_{21} =v_{34} =v_{43} = \frac{1}{1-y} ,\\
     v_{13} =v_{24} =v_{31} =v_{42} =1- \frac{1}{y}\\
     v_{23} =v_{14} =v_{32} =v_{41} =y.
   \end{array} \]
Furthermore, we get
\[ L= \frac{1}{L^{\ast}} = \frac{xy}{x-1} ,M= \frac{1}{M^{\ast}} =
   \frac{(1-y)^{2}}{y^{4}} . \]
It appears that all hyperbolic solutions lie in the first component.

Eliminating $x$ and $y$ between $L$ and $M$ with the condition $x^{2} y^{2} =
(1-x)  (1-y)$ gives the condition (called $A$-polynomial):

{\small{$-L^{8} M+2 \hspace{0.17em} L^{7} M+3 \hspace{0.17em} L^{6} M-2
\hspace{0.17em} L^{5} M+L^{4} M^{2} -6 \hspace{0.17em} L^{4} M+L^{4} -2
\hspace{0.17em} L^{3} M+3 \hspace{0.17em} L^{2} M+2 \hspace{0.17em} LM-M=0$}}

Note that this $A$-polynomial is not the classical one for the group
$\tmop{SL} (2,\mathbb{C})$ but is the same as the one independently found by
C. Zickert {\cite{Zickert-A-pol}}.

Note that hyperbolic and CR-solutions in the first component intersect in a
1-dimensional variety, namely the real part of the hyperbolic solutions. Those
are, from a geometric point of view, degenerate as the group $\tmop{PGL}
(2,\mathbb{R})$ is not linked to a 3-dimensional geometry.

%% file: character-3_mbt.tex
\section{Deforming representations: another way to parametrize the components}\label{deformation}

Having found a specific representation $\rho: \Gamma_8 \to \mathrm{SL(3,\mathbb C)}$ by
algebraic means, there is a complementary method we can often use to determine
the component $\mathcal{V}$ of the variety containing $\rho$. We apply this procedure to a representation in each of the three components found in the previous section. As a product, we get an explicit parametrisation of the representations, i.e. an explicit parametrisation of matrices generating the representation, at least for a Zariski-open subset of the studied component. This information is of course very valuable for the study of the representations found, as in \cite{DF}.

The procedure may be broken
down into four steps, which we now describe briefly; a detailed account of the
method is given in \cite{CLTb}.  Repeated use is made of the LLL algorithm \cite{LLL}
for finding an algebraic number of low degree, close to a number given numerically to
high precision.  For ease of exposition, we shall assume that
$\mathcal{V}$ has dimension $2$; in advance we might not know this dimension,
but usually with experimentation it quickly becomes evident.  Also, the dimension of
the Zariski tangent space at $\rho$ (\cite{CLTa}) provides an upper bound for the
dimension of $\mathcal{V}$.

\begin{itemize}
\item[\bf (i)]  First, we perturb the images of the group generators under $\rho$
slightly, and use these perturbed matrices as a starting point for Newton's method
in conjunction with the group relators, so as to converge numerically to a
representation $\sigma$.  If the original representation $\rho$ is not isolated,
$\sigma$ will typically be some generic representation close to $\rho$; however,
by carefully imposing extra constraints one can ``steer'' the Newton process so
that $\sigma$ has some specified character.  In this way we can obtain a
finite array of representations $\{\sigma_{ij}\}$ such that for each generator $g_k$ of $\Gamma_8$
the traces of the matrices $\alpha_{ij} = \sigma_{ij}(g_k) $ form a suitable array of (numerical
approximations to) rational numbers.  Of course the matrix $\alpha_{ij}$ depends on the
index $k$ of $g_k$, but to avoid excess of notation we omit mention of this index.
\item[\bf (ii)]  Matrices $\beta_{ij}$ are chosen (independently of $k$) so that the conjugates
$\alpha'_{ij} = \beta^{-1}_{ij} \alpha_{ij} \beta_{ij}$ have entries that are all algebraic.
This is always possible, owing to the fact that $\mathcal{V}$ is an algebraic
set; in practice we desire that the degrees of these algebraic numbers should be
as small as possible.
\item[\bf (iii)]  So far, for each group generator $g_k$ we have a two-dimensional array
of matrices $\alpha'_{ij}$ whose entries
are numerical approximations to elements of a number field $\mathcal{F}_{ij}$.
We choose a suitable basis for $\mathcal{F}_{ij}$ over the field of rationals,
noting that the field depends on the parameters $i, j$.  Next, LLL allows us to emerge
from the world of numerical approximations to that of exact quantities by guessing
an expression for each matrix entry as a linear combination of basis elements.
Of course we wish also to be freed from the discrete coordinates $(i, j)$ to
general (complex) coordinates $(u, v)$ for the variety $\mathcal V$, and this is
achieved, up to judicious guesswork, by means of polynomial interpolation across
the $(i,j)$--lattice.

\item[\bf (iv)]  The reader will note that guesswork has been employed at several
points of this procedure, and that substantial use was made of numerical approximations
to algebraic numbers.  However, we have arrived at a function $\Phi$ that associates to
each generator $g_k$ of $\Gamma_8$ an exact matrix, whose entries lie in an algebraic extension
of $\mathbb Q(u,v)$, where $u,v$ are independent transcendentals.  Using software
such as Mathematica or Maple, it is relatively straightforward to verify by formal
manipulation that the matrices $\Phi(g_k)$ satisfy the group relations, {\it i.e.} that
$\Phi$ extends to a representation of $\Gamma_8$.  Following the nomenclature of M. Culler and
P. Shalen, we call $\Phi$ the {\it tautological represention} of $\Gamma_8$ associated to $\mathcal{V}$.
A tautological representation depends of course on a chosen parametrization, {\it i.e.}
a coordinate system for the variety; also, two matrix representations are considered
to be equivalent if they agree up to postmultiplication by an inner automorphism of
the target linear group.
\end{itemize}

The method just described would appear at first sight to be remarkably {\it ad hoc;}\/
however, it often works well.  It was applied to three representations
$\rho_1 \,,\, \rho_2 \,,\, \rho_3$ given in \S \ref{section:example-unipotent} -- one in each of the components found in the previous section. It produces tautological
representations $R_i \; (i=1,2,3)$ of the figure-eight knot group
$\Gamma_8$ into $\mathrm{SL(3,\mathcal{F})}$, where the field $\mathcal{F}$ is an
algebraic extension of $\mathbb Q(u,v)$.  Our brief descriptions of the salient
properties of the $R_i$ are given mostly without proof; however, the reader can
check most of them with the aid of suitable computer algebra software.
\vvv\noind

\subsection{The first component, $R_1$}  \label{ss:R1}
We use the presentation
\[ \Gamma_8 \h=\h \langle a \bc b \;|\; a\,b\,b\,b\,a\,B\,A\,A\,B \rangle \h,\]
where uppercase denotes inverse (see section \ref{section:intro-eight}). The group $\Gamma_8$ is generated by meridians $a\,b \bc a\,b\,b$,
and we may take
\[ m = a\,b \h,\h \ell = a\,b\,A\,B\,a\,B\,A\,b \]
as commuting (meridian, longitude) pair.
Our choice of parameters for the variety is
\[ u = \tr(A) \quad,\quad v = \tr(b) \h,\]
this choice being justified by the fact that
all matrix entries for the tautological representation lie in $\mathbb Q(u,v)$ itself:

\[
a \;\mapsto\; \left[\begin{array}{ccc}0&0&1\\1&0&-u\\0&1&\frac{(-1+v)^2(1+v)}{u}\end{array}\right] \textrm{ and }\]
\[b\;\mapsto\; \left[\begin{array}{ccc}1&
\frac{-u^3(-2+v)+(-1+v)^3(1+v)}{u(-1+v)^2v}&
-\frac{u(-2+v)}{-1+v}\\
0&v&
\frac{-v+v^3}{u}\\
0&
\frac{u}{v-v^2}&
-1\end{array}\right] \quad. \]

For notational convenience, let us identify $a \bc b \bc A \bc B$ with their images
under $R_1$.  These elements have the following traces:
\[ \tr(a) = \frac{(-1+v)^2(1+v)}{u} \h,\h \tr(A) = u \h,\h tr(b) = v \h,\h tr(B) = v \h.\]
On can also compute that
\[ \tr(a\,b) \;=\; \frac{u^3 + (-1 + v)^3(1 + v)^2}{uv(-1 + v)^2} \h,\]
from which it is evident that attempting to assign certain values to the traces
of $A \bc b$ results in singularities.
\vv

The classical hyperbolic representations are those satisfying the
constraint $\tr(a) = \tr(A)$, and one can compute that there are two {\it boundary unipotent}
hyperbolic representations, {\it i.e.} representations mapping peripheral elements
to unipotent matrices.  These are as follows:
\[ (u , v) \;=\; \left(-3 + 2\sqrt{3}\,i \h,\h \frac{1}{2} + \frac{3}{2}\sqrt{3}\,i\right) \quad\textrm{and}\]
\[   (u , v) \;=\; \left(-3 - 2\sqrt{3}\,i \h,\h \frac{1}{2} - \frac{3}{2}\sqrt{3}\,i\right) \h.\]
These form a complex conjugate pair, corresponding to the complete hyperbolic structure
on the figure-eight knot complement.  The deformation space of the complete structure
has one complex dimension, as in the classical $\mathrm{PSL(2,\mathbb C)}$ case.

\subsection{The other two components, $R_2$ and $R_3$}\label{ss:R_2-R_3}

For the other two varieties, we choose the standard ``parabolic'' presentation
\[ \Gamma_8 \h=\h \langle a \bc b \;|\; A\,b\,a\,B\,a\,b\,A\,B\,a\,B \rangle \h,\]
and parameters
\[ u = \tr(A) = \tr(B) \quad,\quad v = \tr(a) = \tr(b) \h.\]
Some matrix entries
lie in a quadratic extension of $\mathbb Q(u,v)$, namely the field generated by
$\sqrt{\Delta}$, where
\[ \Delta \;=\; 4 u^3 + 4 v^3 - u^2 v^2 - 16 u v + 16 \h.\]
Either choice of square root of $\Delta$ yields a valid representation.
We shall see shortly that these varieties are related via precomposition with an
outer automorphism of the knot group.  Here are the tautological representations:
\vvv

\fbox{$R_2$}
\vskip-4ex

{\bfi{12}{21}\rm
\begin{align*}
   a \;\mapsto\; &\left[\begin{array}{ccc}
    \frac{v}{2} & 1 & -\frac{(1-i)(-16 + 8uv - 2v^3 - 4\sqrt{\Delta})}{8u^2 - 6uv^2 + v^4}\\
    \frac{1}{8}(1+i)(-2u+v^2) & \frac{1}{4}(1+i)v & 1\\
    \frac{1}{16}(8 - 4uv + v^3 - 2\sqrt{\Delta}) & \frac{1}{8}(-4u + v^2) & \frac{1}{4}(1-i)v
    \end{array}\right] \\[3ex]
   b\;\mapsto\; &\left[\begin{array}{ccc}
    \frac{v}{2} & i & \frac{(1+i)(-16 + 8uv - 2v^3 - 4\sqrt{\Delta})}{8u^2 - 6uv^2 + v^4}\\
    -\frac{1}{8}(1+i)(-2u+v^2) & \frac{1}{4}(1-i)v & i\\
    -\frac{1}{16}(8 - 4uv + v^3 - 2\sqrt{\Delta}) & -\frac{i}{8}(-4u + v^2) & \frac{1}{4}(1+i)v
    \end{array}\right]
\end{align*}
}

\fbox{$R_3$}
\vskip-4ex

{\bfi{12}{21}\rm
\begin{align*}
   a \;\mapsto\; &\left[\begin{array}{ccc}
    \frac{v}{2} & 1 & -\frac{(1-i)(-16 + 8uv - 2v^3 + 4\sqrt{\Delta})}{8u^2 - 6uv^2 + v^4}\\
    \frac{1}{8}(1+i)(-2u+v^2) & \frac{1}{4}(1+i)v & 1\\
    \frac{1}{16}(8 - 4uv + v^3 + 2\sqrt{\Delta}) & \frac{1}{8}(-4u + v^2) & \frac{1}{4}(1-i)v
    \end{array}\right] \\[3ex]
   b\;\mapsto\; &\left[\begin{array}{ccc}
    \frac{v}{2} & 1 & -\frac{(1-i)(-16 + 8uv - 2v^3 + 4\sqrt{\Delta})}{8u^2 - 6uv^2 + v^4}\\
    \frac{1}{8}(1-i)(-2u+v^2) & \frac{1}{4}(1-i)v & -i\\
    \frac{i}{16}(8 - 4uv + v^3 + 2\sqrt{\Delta}) & \frac{i}{8}(-4u + v^2) & \frac{1}{4}(1+i)v
    \end{array}\right]
\end{align*}
}\ind
\vv

Our first observation is that $R_2$ and $R_3$ satisfy the same relationships with
regard to automorphisms of the knot group $\Gamma_8$ as those set out in Proposition 9.2
of \cite{DF}.  In particular, if $\tau: \Gamma_8 \to \Gamma_8$ is the (non-inner) automorphism
determined by the assignments $a \mapsto a^{-1}\,b\,a \;,\; b \mapsto b$, then
$R_3$ agrees with $R_2 \circ \tau$ up to conjugation by a matrix in
$\mathrm{GL(3,\mathcal F)}$.
\vv

Secondly, the image of $R_2$ (also that of $R_3$) is a copy of the $(3,3,4)$--triangle group
\[ T(3,3,4) \;=\; \langle \alpha \bc \beta \bc \gamma \;|\;
   \alpha^3 = \beta^3 = \gamma^4 = \alpha\beta\gamma = 1 \rangle \h.\]
Specifically, if $c$ denotes the commutator $a\,b\,A\,B$\,, we have
\[ R_2(c) \;=\; \left[\begin{array}{ccc}1&0&0\\0&i&0\\0&0&-i\end{array}\right] \h,\]
independently of $u \bc v$, from which it follows immediately that $R_2(c)$ has order $4$.
Moreover, one can verify that the images under $R_2$ of $a\,c \;,\; c\,a\,c$ each have
order $3$.  Since $a \bc c$ generate the knot group $\Gamma_8$, the result follows by taking
$\alpha = R_2(a\,c) \;,\; \beta = R_2(c\,a\,c)^{-1} \;,\; \gamma = R_2(c)$.  It was
previously known that the images of the representations $\rho_2 \bc \rho_3$
were isomorphic to $T(3,3,4)$, see \cite{DF}.

\subsection{A focus on the $\mathrm{PU}(2,1)$ case}

As studied e.g. in \cite{DF}, the representations with image inside the real group $\mathrm{PU}(2,1)$ may carry interesting geometric properties. So we seek representations preserving Hermitian forms.

In the case of the first component $R_1$, the following condition on traces holds:
\[ \tr(A) \,=\, \overline{\tr(a)} \quad,\quad \tr(B) \,=\, \overline{\tr(b)} \h.\]
In particular, we note that since $\tr(B) = \tr(b) = v$, the parameter $v$ is forced to be real
for such representations.  Then, from $\tr(A) = \tr(a)$, we obtain
\[ 
\frac{(-1+v)^2(1+v)}{\overline{u}} \;=\; u\;,
\quad\text{or}\quad
u \;=\; |  \left(-1+v\right)\sqrt{1+v} |  \, z ,
\]
with $z$ lying on the unit circle.
\vv

It is in fact possible to compute a Hermitian matrix $H$ as a function of the coordinates
$v \bc z$,
\[ H(v,z) \;=\; \left[\begin{array}{ccc}1&r&s\\\overline{r}&1&r\\\overline{s}&\overline{r}&1\end{array}\right]\h,\]
whose associated form is respected by the representations just described;
the entries of $H$ lie in an Abelian extension of degree $4$ of $\mathbb Q(v,\mathrm{Re}(z))$,
and are a little too cumbersome for inclusion here.  However, for specific $v \bc z$
the signature of $H$ enables one to determine whether the representation
takes values in $\mathrm{SU(3)}$ or $\mathrm{SU(2,1)}$, or whether the associated
form is degenerate.  Indeed, the representations with degenerate forms occupy some
isolated points in the $(v,z)$--variety, together with some simple closed curves that
separate $\mathrm{SU(3)}$--representations from $\mathrm{SU(2,1)}$--representations.
\vv

One can compute that the variety $R_1$ gives rise to a single boundary
unipotent representation into $\mathrm{PU(2,1)}$, obtained by composing the representation
\[ (u , v) \;=\; (-\sqrt{3}\,i \;,\; 2) \]
with the natural projectivization map $\mathrm{SU(2,1) \to PU(2,1)}$.  The variety
$R_1$ also contains the other two lifts of this $\mathrm{PU(2,1)}$--
representation to $\mathrm{SU(2,1)}$, at
\[ (u , v) \;=\; (-\omega(\sqrt{3}\,i) \;,\; 2) \, ,\quad
   (u , v) \;=\; (-\omega^2(\sqrt{3}\,i) \;,\; 2) \quad
   (\text{where}\quad \omega = e^{2\pi i/3}).\]

\bigskip

The other two components $R_2$ and $R_3$ each give rise to just one boundary unipotent representation into
$\mathrm{PU(2,1)}$, given by $(u,v) = (3,3)$.
\vv

A necessary condition on traces for a representation in $R_2$ or $R_3$ to respect a Hermitian
forms is
\[ \tr(a) = \tr(A) \;,\; \tr(b) = \tr(B) \quad,\quad\text{\it i.e.}\quad v = \overline{u} \h.\]

Using the same technique as that described at the beginning of this section,
{\it i.e.} invoking LLL together with polynomial interpolation, the following diagonal matrix:
\[ H \;=\; \left[\begin{array}{ccc}
-\frac{1}{8}(\Delta - 16)\left(\sqrt{\Delta} + |u|^2 - 4 \right) & 0 & 0\\
0 & \Delta - 16 & 0\\
0 & 0 & -8\left(\sqrt{\Delta} + 4 \right)
\end{array}\right] \]
was found, satisfying
\[ H a^{-1} \;=\; a^{\,\dagger}H \h,\h H b^{-1} \;=\; b^{\,\dagger}H
\quad (v = \overline{u} \;,\; \Delta \geq 0) \h,\]
where the superscript $^\dagger$ denotes conjugate transpose, and
$\Delta = 4 u^3 + 4 v^3 - u^2 v^2 - 16 u v + 16 \;=\;
4(u^3 + \overline{u}^3) - |u|^4 - 16|u|^2 + 16$.
Clearly the condition $\Delta \geq 0$ is required for $H$ to be Hermitian.  We note that
the form represented by $H$ degenerates at $\Delta = 16$; in the present context $v = \overline{u}$
the corresponding values of $u$ are $0 \,,\, 4 \,,\, 4e^{2\pi i/3} \,,\, 4e^{4\pi i/3}$.
\vv

Substituting $u = x + iy \bc v = x - iy \h (x \bc y \in \mathbb R)$ in the expression for
$\Delta(u,v)$, we obtain (with mild abuse of notation)
\[ \Delta(x,y) \;=\; -x^4 - y^4 - 2x^2 y^2 - 24x y^2 + 8x^3 - 16x^2 - 16y^2 + 16 \h.\]
\ind

We have therefore established that
representations respecting Hermitian forms occur in the bounded complementary regions
in the $(x,y)$--plane of the curve $\Delta(x,y) = 0$, see Figure \ref{fig:su(2,1)}.
The $\mathrm{SU(3)}$ representations lie in the central complementary region, and
the $\mathrm{SU(2,1)}$ representations lie in the three petals.  Corresponding
Hermitian representations occur also on the curve itself, except for singularities
at the three self-intersection points of the curve, where representations are
not defined.  The quantity $\Delta$ has
geometric significance in that it is related to the Brehm shape invariant \cite{Brehm}
of the triangle in the associated homogeneous space (complex projective plane
$\mathbb{CP}^2$ or complex hyperbolic plane $\mathbb{CH}^2$), whose vertices are the
fixed points of the elliptic isometries $\alpha \bc \beta \bc \gamma$.

\begin{figure}[ht]
\includegraphics[scale=1]{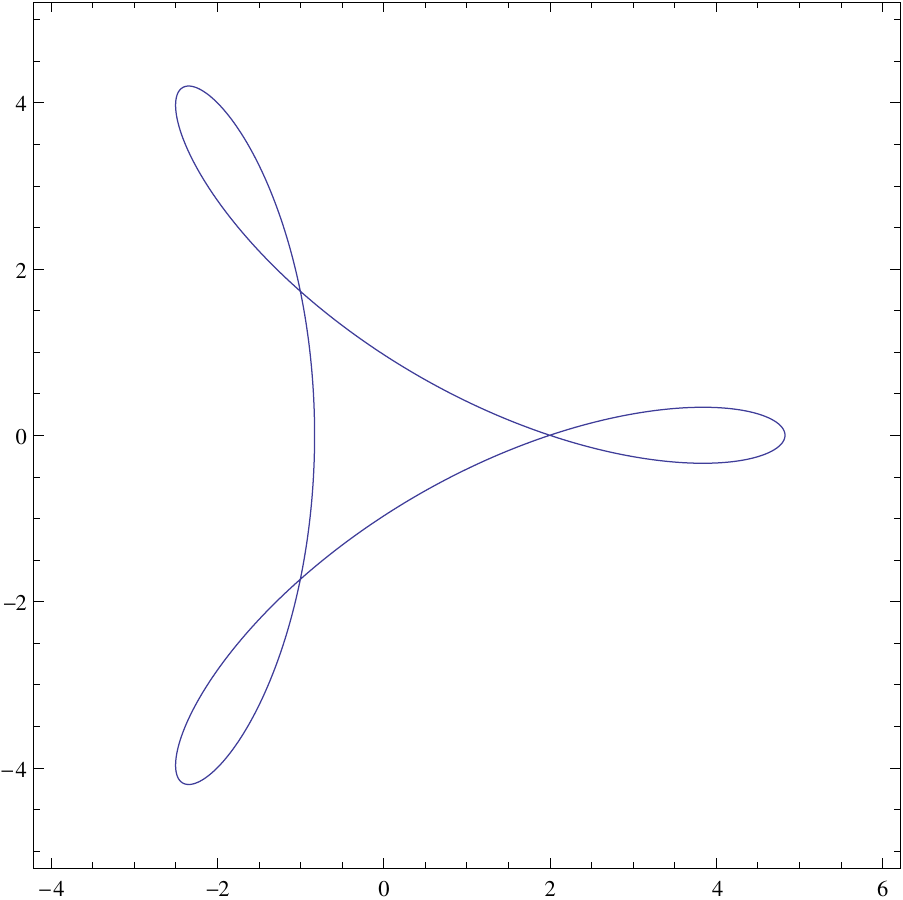}
\caption{$\mathrm{SU(3)}$ and $\mathrm{SU(2,1)}$ representations in $R_2$}\label{fig:su(2,1)}
\end{figure}

\bibliography{bibli-1}
\bibliographystyle{plain}